# Linear Mappings of Quaternion Algebra

Aleks Kleyn

ABSTRACT. In the paper I considered linear and antilinear automorphisms of quaternion algebra. I proved the theorem that there is unique expansion of $R$-linear mapping of quaternion algebra relative to the given set of linear and antilinear automorphisms.

## Contents

## 1. Preface to Version 1

This paper by its form looks like reference book for elementary functions. However, the highlight of the paper is the theorem 12.1 which is a powerful source of new ideas.

When I was writing the paper [2], I hoped until the last minute that I would be able to prove that $D$-linear mapping in $D$-algebra $A$ with conjugation has expansion into sum of $A$-linear and $A$-antilinear mappings. Since it is necessary to write down the system of linear equations to solve this problem, then I finally realized that even in case of quaternion algebra this problem has a negative answer.

Aleks_Kleyn@MailAPS.org.
http://sites.google.com/site/AleksKleyn/.
http://arxiv.org/a/kleyn_a_1.
http://AleksKleyn.blogspot.com/.

The statement that there exists nontrivial linear automorphism in quaternion algebra was a complete surprise to me.[1] So it took me a while to understand how this statement is important.

When I started research of algebra with a countable basis, suddenly I realised that if I have 3 linear automorphisms (namely, $\overline{E}$, $\overline{E}_1$, $\overline{E}_2$ ) and antilinear automorphism $\overline{I}$, then there is hope that for any $R$-linear mapping[2] $f$ there exists unique expansion

$$(1.1) \qquad f = a_0 \circ \overline{E} + a_1 \circ \overline{E}_1 + a_2 \circ \overline{E}_2 + a_3 \circ \overline{I}$$

The solution seemed simple, however corresponding system of linear equations is singular. Various attempts to modify the third term were unsuccessful.

Fortune smiled me when I asked myself whether the mapping $\overline{E}_3$ is linear automorphism.

## 2. Preface to Version 2

Anthony Sudbery drew my attention to the statement that any transformation of the form

$$(2.1) \qquad x \to a\,x\,a^{-1} \quad a \in H$$

is linear automorphism.[3] In particular, the mapping $\overline{E}_1$ has form

$$\overline{E}_1 \circ x = (1 + \sqrt{3}(i+j+k))x(1 + \sqrt{3}(i+j+k))^{-1}$$

Can we prove that any linear automorphism has form (2.1)?

The theorem 12.1 has another and very important consequence. In the section [1]-2.4, I considered $D$-vector space (in this case, $D = H$). It follows from the theorem 12.1 that dimension of $H$-vector space $V$ does not coincide with dimension of $H\star$-vector space $V$.

## 3. Conventions

**Convention 3.1.** *Let $A$ be free finite dimensional algebra. Considering expansion of element of algebra $A$ relative basis $\overline{\overline{e}}$ we use the same root letter to denote this element and its coordinates. However we do not use vector notation in algebra. In expression $a^2$, it is not clear whether this is component of expansion of element $a$ relative basis, or this is operation $a^2 = aa$. To make text clearer we use separate color for index of element of algebra. For instance,*

$$a = a^i \overline{e}_i$$

□

**Convention 3.2.** *If free finite dimensional algebra has unit, then we identify the vector of basis $\overline{e}_0$ with unit of algebra.* □

---

[1]Sometimes it is very hard to reproduce the actual sequence of events. However, when I was preparing an example [2]-7.2, I recalled the diagram of symmetry of product in quaternion algebra; the diagram was given in [4], p. 7.

[2]Since $R$-linear mappings are main tool in my research of calculus, I usually use the term linear mapping. However in case when we consider several algebras at the same time, it is necessary to specify algebra over which the mapping is linear.

[3]Such mappings are well known (see, for instance, section [3]-3.2, as well p. [5]-643). However I did not think about this case because my attention was focused on mappings of the form $y = af(x)$.

**Convention 3.3.** *For given D-algebra A we define left shift* $a\circ$ *by the equation*

$$a\circ x = ax$$

*and right shift* $a\star$ *by the equation*

$$a\star x = xa$$

□

Without a doubt, the reader may have questions, comments, objections. I will appreciate any response.

## 4. Linear Automorphism of Quaternion Algebra

**Theorem 4.1.** *Coordinates of linear automorphism of quaternion algebra satisfy to the system of equations*

(4.1)
$$\begin{cases} r_1^1 = r_2^2 r_3^3 - r_2^3 r_3^2 & r_2^1 = r_3^2 r_1^3 - r_3^3 r_1^2 & r_3^1 = r_1^2 r_2^3 - r_1^3 r_2^2 \\ r_1^2 = r_2^3 r_3^1 - r_2^1 r_3^3 & r_2^2 = r_3^3 r_1^1 - r_3^1 r_1^3 & r_3^2 = r_1^3 r_2^1 - r_1^1 r_2^3 \\ r_1^3 = r_2^1 r_3^2 - r_2^2 r_3^1 & r_2^3 = r_3^1 r_1^2 - r_3^2 r_1^1 & r_3^3 = r_1^1 r_2^2 - r_1^2 r_2^1 \end{cases}$$

*Proof.* According to the theorems [3]-3.3.1, [2]-6.4, linear automorphism of quaternion algebra satisfies to equations

(4.2)
$$\begin{aligned} r_0^l &= r_0^p r_0^q C_{pq}^l & r_1^l &= r_0^p r_1^q C_{pq}^l & r_2^l &= r_0^p r_2^q C_{pq}^l & r_3^l &= r_0^p r_3^q C_{pq}^l \\ r_1^l &= r_1^p r_0^q C_{pq}^l & -r_0^l &= r_1^p r_1^q C_{pq}^l & r_3^l &= r_1^p r_2^q C_{pq}^l & -r_2^l &= r_1^p r_3^q C_{pq}^l \\ r_2^l &= r_2^p r_0^q C_{pq}^l & -r_3^l &= r_2^p r_1^q C_{pq}^l & -r_0^l &= r_2^p r_2^q C_{pq}^l & r_0^l &= r_2^p r_3^q C_{pq}^l \\ r_3^l &= r_3^p r_0^q C_{pq}^l & r_2^l &= r_3^p r_1^q C_{pq}^l & -r_1^l &= r_3^p r_2^q C_{pq}^l & -r_0^l &= r_3^p r_3^q C_{pq}^l \end{aligned}$$

From the equation (4.2), it follows that

(4.3)
$$\begin{aligned} r_1^l &= r_0^p r_1^q C_{pq}^l = r_2^p r_3^q C_{pq}^l & r_0^p r_1^q C_{pq}^l &= r_0^p r_1^q C_{qp}^l & r_2^p r_3^q C_{pq}^l &= -r_2^p r_3^q C_{qp}^l \\ r_2^l &= r_0^p r_2^q C_{pq}^l = r_3^p r_1^q C_{pq}^l & r_0^p r_2^q C_{pq}^l &= r_0^p r_2^q C_{qp}^l & r_3^p r_1^q C_{pq}^l &= -r_1^p r_3^q C_{pq}^l \\ r_3^l &= r_0^p r_3^q C_{pq}^l = r_1^p r_2^q C_{pq}^l & r_0^p r_3^q C_{pq}^l &= r_0^p r_3^q C_{qp}^l & r_1^p r_2^q C_{pq}^l &= -r_1^p r_2^q C_{qp}^l \end{aligned}$$

(4.4) $$r_0^l = r_0^p r_0^q C_{pq}^l = -r_1^p r_1^q C_{pq}^l = -r_2^p r_2^q C_{pq}^l = -r_3^p r_3^q C_{pq}^l$$

If $l = 0$, then from the equation

$$C_{pq}^0 = C_{qp}^0$$

it follows that

(4.5) $$r_i^p r_j^q C_{pq}^0 = r_i^p r_j^q C_{qp}^0$$

From the equation (4.3) for $l = 0$ and the equation (4.5), it follows that

(4.6) $$r_1^0 = r_2^0 = r_3^0 = 0$$

If $l = 1, 2, 3,$ then we can write the equation (4.3) in the following form

$$(4.7) \begin{cases} r_i^l = r_0^l r_i^0 C_{l0}^l + r_0^0 r_i^l C_{0l}^l + r_0^a r_i^b C_{ab}^l + r_0^b r_i^a C_{ba}^l \\ \quad r_0^l r_i^0 C_{l0}^l + r_0^0 r_i^l C_{0l}^l + r_0^a r_i^b C_{ab}^l + r_0^b r_i^a C_{ba}^l \\ = r_0^l r_i^0 C_{0l}^l + r_0^0 r_i^l C_{l0}^l + r_0^a r_i^b C_{ba}^l + r_0^b r_i^a C_{ab}^l \\ \quad i = 1, 2, 3 \\ r_i^l = r_k^0 r_j^l C_{0l}^l + r_k^l r_j^0 C_{l0}^l + r_k^a r_j^b C_{ab}^l + r_k^b r_j^a C_{ba}^l \\ \quad r_k^0 r_j^l C_{0l}^l + r_k^l r_j^0 C_{l0}^l + r_k^a r_j^b C_{ab}^l + r_k^b r_j^a C_{ba}^l \\ = -r_k^0 r_j^l C_{l0}^l - r_k^l r_j^0 C_{0l}^l - r_k^a r_j^b C_{ba}^l - r_k^b r_j^a C_{ab}^l \\ \quad i = 1 \quad k = 2 \quad j = 3 \\ \quad i = 2 \quad k = 3 \quad j = 1 \\ \quad i = 3 \quad k = 1 \quad j = 2 \\ \quad 0 < a < b \quad a \neq l \quad b \neq l \end{cases}$$

From equations (4.7), (4.6) and equations

$$(4.8) \begin{aligned} C_{0l}^l = C_{l0}^l = 1 \\ C_{ab}^l = -C_{ba}^l \end{aligned}$$

it follows that

$$(4.9) \begin{cases} r_i^l = r_0^0 r_i^l + r_0^a r_i^b C_{ab}^l - r_0^b r_i^a C_{ab}^l \\ \quad r_0^0 r_i^l + r_0^a r_i^b C_{ab}^l - r_0^b r_i^a C_{ab}^l \\ = r_0^0 r_i^l - r_0^a r_i^b C_{ab}^l + r_0^b r_i^a C_{ab}^l \\ \quad i = 1, 2, 3 \\ r_i^l = r_k^a r_j^b C_{ab}^l - r_k^b r_j^a C_{ab}^l \\ \quad r_k^a r_j^b C_{ab}^l - r_k^b r_j^a C_{ab}^l \\ = r_k^a r_j^b C_{ab}^l - r_k^b r_j^a C_{ab}^l \\ \quad i = 1 \quad k = 2 \quad j = 3 \\ \quad i = 2 \quad k = 3 \quad j = 1 \\ \quad i = 3 \quad k = 1 \quad j = 2 \\ \quad 0 < a < b \quad a \neq l \quad b \neq l \end{cases}$$

From equations (4.9) it follows that

(4.10)
$$\begin{cases} r_i^l = r_0^0 r_i^l \\ \quad r_0^a r_i^b - r_0^b r_i^a = 0 \\ \quad i = 1, 2, 3 \\ r_i^l = r_k^a r_j^b C_{ab}^l - r_k^b r_j^a C_{ab}^l \\ \quad i = 1 \quad k = 2 \quad j = 3 \\ \quad i = 2 \quad k = 3 \quad j = 1 \\ \quad i = 3 \quad k = 1 \quad j = 2 \\ \quad 0 < a < b \quad a \neq l \quad b \neq l \end{cases}$$

From the equation (4.10) it follows that

(4.11) $$r_0^0 = 1$$

From the equation (4.4) for $l = 0$, it follows that

(4.12)
$$\begin{aligned} r_0^0 &= r_0^0 r_0^0 - r_0^1 r_0^1 - r_0^2 r_0^2 - r_0^3 r_0^3 \\ &= -r_i^0 r_i^0 + r_i^1 r_i^1 + r_i^2 r_i^2 + r_i^3 r_i^3 \\ i &= 1, 2, 3 \end{aligned}$$

From equations (4.6), (4.10), (4.12), it follows that

(4.13)
$$\begin{aligned} 0 &= r_0^1 r_0^1 + r_0^2 r_0^2 + r_0^3 r_0^3 \\ 1 &= r_i^1 r_i^1 + r_i^2 r_i^2 + r_i^3 r_i^3 \\ i &= 1, 2, 3 \end{aligned}$$

From equations (4.13) it follows that[4]

(4.14) $$r_0^1 = r_0^2 = r_0^3 = 0$$

From the equation (4.4) for $l > 0$, it follows that

(4.15)
$$\begin{aligned} r_0^l &= r_0^l r_0^0 C_{l0}^l + r_0^0 r_0^l C_{0l}^l + r_0^a r_0^b C_{ab}^l + r_0^b r_0^a C_{ba}^l \\ &= -r_i^l r_i^0 C_{l0}^l - r_i^0 r_i^l C_{0l}^l - r_i^a r_i^b C_{ab}^l - r_i^b r_i^a C_{ba}^l \\ i &> 0 \\ l &> 0 \quad 0 < a < b \quad a \neq l \quad b \neq l \end{aligned}$$

---

[4]Here, we rely on the fact that quaternion algebra is defined over real field. If we consider quaternion algebra over complex field, then the equation (4.13) defines a cone in the complex space. Correspondingly, we have wider choice of coordinates of linear automorphism.

Equations (4.15) are identically true by equations (4.6), (4.14), (4.8). From equations (4.14), (4.10), it follows that

$$(4.16) \quad \begin{cases} r_i^l = r_k^a r_j^b C_{ab}^l - r_k^b r_j^a C_{ab}^l \\ \quad i = 1 \quad k = 2 \quad j = 3 \\ \quad i = 2 \quad k = 3 \quad j = 1 \\ \quad i = 3 \quad k = 1 \quad j = 2 \\ \quad l > 0 \quad 0 < a < b \quad a \neq l \quad b \neq l \end{cases}$$

Equations (4.1) follow from equations (4.16). □

## 5. Mapping $\overline{E}$

It is evident that coordinates of mapping

$$\overline{E} : H \to H \quad \overline{E} \circ x = x$$

$$E_j^i = \delta_j^i$$

satisfy the equation (4.1).

**Theorem 5.1.** *We can identify the mapping*

$$a \circ \overline{E} : H \to H \quad a \in H$$

*and matrix*

$$(5.1) \quad J_{l \cdot a} = \begin{pmatrix} a^0 & -a^1 & -a^2 & -a^3 \\ a^1 & a^0 & -a^3 & a^2 \\ a^2 & a^3 & a^0 & -a^1 \\ a^3 & -a^2 & a^1 & a^0 \end{pmatrix}$$

*Proof.* The product of quaternions

$$a = a^0 + a^1 i + a^2 j + a^3 k$$

and

$$x = x^0 + x^1 i + x^2 j + x^3 k$$

has form

$$ax = a^0 x^0 - a^1 x^1 - a^2 x^2 - a^3 x^3 + (a^0 x^1 + a^1 x^0 + a^2 x^3 - a^3 x^2)i$$
$$+ (a^0 x^2 + a^2 x^0 + a^3 x^1 - a^1 x^3)j + (a^0 x^3 + a^3 x^0 + a^1 x^2 - a^2 x^1)k$$

Therefore, function $f_a(x) = ax$ has Jacobian matrix (5.1). □

**Theorem 5.2.** *We can identify the mapping*

$$a \star \overline{E} : H \to H \quad a \in H$$

*and matrix*

$$(5.2) \quad J_{r \cdot a} = \begin{pmatrix} a^0 & -a^1 & -a^2 & -a^3 \\ a^1 & a^0 & a^3 & -a^2 \\ a^2 & -a^3 & a^0 & a^1 \\ a^3 & a^2 & -a^1 & a^0 \end{pmatrix}$$

*Proof.* The product of quaternions

$$x = x^0 + x^1 i + x^2 j + x^3 k$$

and

$$a = a^0 + a^1 i + a^2 j + a^3 k$$

has form

$$xa = x^0 a^0 - x^1 a^1 - x^2 a^2 - x^3 a^3 + (x^0 a^1 + x^1 a^0 + x^2 a^3 - x^3 a^2)i$$
$$+ (x^0 a^2 + x^2 a^0 + x^3 a^1 - x^1 a^3)j + (x^0 a^3 + x^3 a^0 + x^1 a^2 - x^2 a^1)k$$

Therefore, function $f_a(x) = ax$ has Jacobian matrix (5.2). □

## 6. Mapping $\overline{E}_1$

**Theorem 6.1.** *Coordinates of mapping*

$$\overline{E}_1 : H \to H \quad \overline{E}_1 \circ x = x^0 + x^2 i + x^3 j + x^1 k$$

$$E_{1 \cdot 0}^{0} = 1 \quad E_{1 \cdot 2}^{1} = 1 \quad E_{1 \cdot 3}^{2} = 1 \quad E_{1 \cdot 1}^{3} = 1$$

$$E_1 = \begin{pmatrix} 1 & 0 & 0 & 0 \\ 0 & 0 & 1 & 0 \\ 0 & 0 & 0 & 1 \\ 0 & 1 & 0 & 0 \end{pmatrix}$$

*satisfy the equation* (4.1).

*Proof.* We can verify directly that (4.1) are true. However, it is easy to verify that using of permutation

$$\begin{pmatrix} 1 & 2 & 3 \\ 2 & 3 & 1 \end{pmatrix}$$

over the set of lower indices preserves the set of equations. □

*Remark* 6.2. It can be verified directly that $\overline{E}_1$ is linear automorphism. Let

$$a = a^0 + a^1 i + a^2 j + a^3 k$$
$$b = b^0 + b^1 i + b^2 j + b^3 k$$

Then

$$ab = a^0 b^0 - a^1 b^1 - a^2 b^2 - a^3 b^3 + (a^0 b^1 + a^1 b^0 + a^2 b^3 - a^3 b^2)i$$
$$+ (a^0 b^2 + a^2 b^0 + a^3 b^1 - a^1 b^3)j + (a^0 b^3 + a^3 b^0 + a^1 b^2 - a^2 b^1)k$$
$$\overline{E}_3 \circ a = a^0 + a^2 i + a^1 j - a^3 k$$

$$\overline{E}_3 \circ b = b^0 + b^2 i + b^1 j - b^3 k$$

$$(\overline{E}_3 \circ a)(\overline{E}_3 \circ b) = a^0 b^0 - a^2 b^2 - a^1 b^1 - a^3 b^3 + (a^0 b^2 + a^2 b^0 - a^1 b^3 + a^3 b^1) i$$
$$+ (a^0 b^1 + a^1 b^0 - a^3 b^2 + a^2 b^3) j + (-a^0 b^3 - a^3 b^0 + a^2 b^1 - a^1 b^2) k$$
$$= \overline{E}_3 \circ (ab)$$

□

**Theorem 6.3.** *The mapping $\overline{E}_1$ has form*

$$\overline{E}_1 \circ a = \frac{1}{4}(a - iai - jaj - kak - ia + ai - kaj - jak$$
$$-ja - kai + aj - iak - ka - jai - iaj + ak)$$

*Proof.* According to the theorem [3]-3.3.4, standard components of the mapping $\overline{E}_1$ have form

$$\begin{array}{cccc}
E_1^{00} = \frac{1}{4} & E_1^{11} = -\frac{1}{4} & E_1^{22} = -\frac{1}{4} & E_1^{33} = -\frac{1}{4} \\
E_1^{10} = -\frac{1}{4} & E_1^{01} = \frac{1}{4} & E_1^{32} = -\frac{1}{4} & E_1^{23} = -\frac{1}{4} \\
E_1^{20} = -\frac{1}{4} & E_1^{31} = -\frac{1}{4} & E_1^{02} = \frac{1}{4} & E_1^{13} = -\frac{1}{4} \\
E_1^{30} = -\frac{1}{4} & E_1^{21} = -\frac{1}{4} & E_1^{12} = -\frac{1}{4} & E_1^{03} = \frac{1}{4}
\end{array}$$

□

**Theorem 6.4.** *We can identify the mapping*

$$a \circ \overline{E}_1 : H \to H \quad a \in H$$

*and matrix*

(6.1) $$J_{1l \cdot a} = \begin{pmatrix} a^0 & -a^3 & -a^1 & -a^2 \\ a^1 & a^2 & a^0 & -a^3 \\ a^2 & -a^1 & a^3 & a^0 \\ a^3 & a^0 & -a^2 & a^1 \end{pmatrix}$$

(6.2) $$J_{1l \cdot a} = J_{l \cdot a} \circ E_1$$

*Proof.* The product of quaternions

$$a = a^0 + a^1 i + a^2 j + a^3 k$$

and

$$\overline{E}_1 \circ x = x^0 + x^2 i + x^3 j + x^1 k$$

has form

$$a \circ \overline{E}_1 \circ x = a^0 x^0 - a^1 x^2 - a^2 x^3 - a^3 x^1 + (a^0 x^2 + a^1 x^0 + a^2 x^1 - a^3 x^3) i$$
$$+ (a^0 x^3 + a^2 x^0 + a^3 x^2 - a^1 x^1) j + (a^0 x^1 + a^3 x^0 + a^1 x^3 - a^2 x^2) k$$

Therefore, function $f_a \circ x = a \circ \overline{E}_1 \circ x$ has Jacobian matrix (6.1). The equation (6.2) follows from the chain of equations

$$J_{l \cdot a \circ} \circ E_2 = \begin{pmatrix} a^0 & -a^1 & -a^2 & -a^3 \\ a^1 & a^0 & -a^3 & a^2 \\ a^2 & a^3 & a^0 & -a^1 \\ a^3 & -a^2 & a^1 & a^0 \end{pmatrix} \circ \begin{pmatrix} 1 & 0 & 0 & 0 \\ 0 & 0 & 1 & 0 \\ 0 & 0 & 0 & 1 \\ 0 & 1 & 0 & 0 \end{pmatrix}$$

$$= \begin{pmatrix} a^0 & -a^3 & -a^1 & -a^2 \\ a^1 & a^2 & a^0 & -a^3 \\ a^2 & -a^1 & a^3 & a^0 \\ a^3 & a^0 & -a^2 & a^1 \end{pmatrix} = J_{1l \cdot a}$$

□

## 7. Mapping $\overline{E}_2$

**Theorem 7.1.** *Since mappings*

$$f_1 : H \to H$$
$$f_2 : H \to H$$

*are linear automorphisms of quaternion algebra, then the mapping* $f_2 \circ f_1$ *is linear automorphism of quaternion algebra.*

*Proof.* The statement of the theorem is corollary of following equations

$$f_2 \circ f_1 \circ (a+b) = f_2 \circ (f_1 \circ a + f_1 \circ b) = f_2 \circ f_1 \circ a + f_2 \circ f_1 \circ b$$
$$f_2 \circ f_1 \circ (ra) = f_2 \circ (r(f_1 \circ a)) = r(f_2 \circ f_1 \circ a)$$
$$f_2 \circ f_1 \circ (ab) = f_2((f_1 \circ a)(f_1 \circ b)) = (f_2 \circ f_1 \circ a)(f_2 \circ f_1 \circ b)$$

□

**Theorem 7.2.** *The mapping*

$$\overline{E}_2 : H \to H \quad \overline{E}_2(x) = x^0 + x^3 i + x^1 j + x^2 k$$

$$E_{2 \cdot 0}^{0} = 1 \quad E_{2 \cdot 3}^{1} = 1 \quad E_{2 \cdot 1}^{2} = 1 \quad E_{2 \cdot 2}^{3} = 1$$

$$E_2 = \begin{pmatrix} 1 & 0 & 0 & 0 \\ 0 & 0 & 0 & 1 \\ 0 & 1 & 0 & 0 \\ 0 & 0 & 1 & 0 \end{pmatrix}$$

*is the linear automorphism of quaternion algebra.*

*Proof.* The statement of the theorem is corollary of the theorem 7.1 and the equation

$$(E_1)^2 = E_2$$

□

**Theorem 7.3.** *The mapping $\overline{E}_2$ has form*

$$\overline{E}_2 \circ a = \frac{1}{4}(a - iai - jaj - kak + ia - ai - kaj - jak$$
$$+ja - kai - aj - iak + ka - jai - iaj - ak)$$

*Proof.* According to the theorem [3]-3.3.4, standard components of the mapping $\overline{E}_2$ have form

$$\begin{array}{llll}
E_2^{00} = \frac{1}{4} & E_2^{11} = -\frac{1}{4} & E_2^{22} = -\frac{1}{4} & E_2^{33} = -\frac{1}{4} \\
E_2^{10} = \frac{1}{4} & E_2^{01} = -\frac{1}{4} & E_2^{32} = -\frac{1}{4} & E_2^{23} = -\frac{1}{4} \\
E_2^{20} = \frac{1}{4} & E_2^{31} = -\frac{1}{4} & E_2^{02} = -\frac{1}{4} & E_2^{13} = -\frac{1}{4} \\
E_2^{30} = \frac{1}{4} & E_2^{21} = -\frac{1}{4} & E_2^{12} = -\frac{1}{4} & E_2^{03} = -\frac{1}{4}
\end{array}$$

□

**Theorem 7.4.** *We can identify the mapping*

$$a \circ \overline{E}_2 : H \to H \quad a \in H$$

*and matrix*

(7.1) $$J_{2l \cdot a} = \begin{pmatrix} a^0 & -a^2 & -a^3 & -a^1 \\ a^1 & -a^3 & a^2 & a^0 \\ a^2 & a^0 & -a^1 & a^3 \\ a^3 & a^1 & a^0 & -a^2 \end{pmatrix}$$

(7.2) $$J_{2l \cdot a} = J_{l \cdot a} \circ E_2$$

*Proof.* The product of quaternions

$$a = a^0 + a^1 i + a^2 j + a^3 k$$

and

$$\overline{E}_2 \circ x = x^0 + x^3 i + x^1 j + x^2 k$$

has form

$$a \circ \overline{E}_2 \circ x = a^0 x^0 - a^1 x^3 - a^2 x^1 - a^3 x^2 + (a^0 x^3 + a^1 x^0 + a^2 x^2 - a^3 x^1)i$$
$$+ (a^0 x^1 + a^2 x^0 + a^3 x^3 - a^1 x^2)j + (a^0 x^2 + a^3 x^0 + a^1 x^1 - a^2 x^3)k$$

Therefore, function $f_a \circ x = a \circ \overline{E}_1 \circ x$ has Jacobian matrix (7.1). The equation (7.2) follows from the chain of equations

$$J_{l \cdot a} \circ E_2 = \begin{pmatrix} a^0 & -a^1 & -a^2 & -a^3 \\ a^1 & a^0 & -a^3 & a^2 \\ a^2 & a^3 & a^0 & -a^1 \\ a^3 & -a^2 & a^1 & a^0 \end{pmatrix} \circ \begin{pmatrix} 1 & 0 & 0 & 0 \\ 0 & 0 & 0 & 1 \\ 0 & 1 & 0 & 0 \\ 0 & 0 & 1 & 0 \end{pmatrix}$$

$$= \begin{pmatrix} a^0 & -a^2 & -a^3 & -a^1 \\ a^1 & -a^3 & a^2 & a^0 \\ a^2 & a^0 & -a^1 & a^3 \\ a^3 & a^1 & a^0 & -a^2 \end{pmatrix} = J_{2l \cdot a}$$

□

## 8. Mapping $\overline{E}_3$

**Theorem 8.1.** *Coordinates of mapping*

$$\overline{E}_3 : H \to H \quad \overline{E}_3 \circ x = x^0 + x^2 i + x^1 j - x^3 k$$

$$E_{3 \cdot 0}^{\ 0} = 1 \quad E_{3 \cdot 2}^{\ 1} = 1 \quad E_{3 \cdot 1}^{\ 2} = 1 \quad E_{3 \cdot 3}^{\ 3} = -1$$

$$E_3 = \begin{pmatrix} 1 & 0 & 0 & 0 \\ 0 & 0 & 1 & 0 \\ 0 & 1 & 0 & 0 \\ 0 & 0 & 0 & -1 \end{pmatrix}$$

*satisfy the equation* (4.1).

*Proof.* We can verify directly that (4.1) are true

$$E_{3 \cdot 1}^{\ 1} = E_{3 \cdot 2}^{\ 2} E_{3 \cdot 3}^{\ 3} - E_{3 \cdot 2}^{\ 3} E_{3 \cdot 3}^{\ 2} = 0 - 0 = 0$$
$$E_{3 \cdot 2}^{\ 1} = E_{3 \cdot 3}^{\ 2} E_{3 \cdot 1}^{\ 3} - E_{3 \cdot 3}^{\ 3} E_{3 \cdot 1}^{\ 2} = 0 - (-1)1 = 1$$
$$E_{3 \cdot 3}^{\ 1} = E_{3 \cdot 1}^{\ 2} E_{3 \cdot 2}^{\ 3} - E_{3 \cdot 1}^{\ 3} E_{3 \cdot 2}^{\ 2} = 1 * 0 - 0 = 0$$
$$E_{3 \cdot 1}^{\ 2} = E_{3 \cdot 2}^{\ 3} E_{3 \cdot 3}^{\ 1} - E_{3 \cdot 2}^{\ 1} E_{3 \cdot 3}^{\ 3} = 0 - 1(-1) = 1$$
$$E_{3 \cdot 2}^{\ 2} = E_{3 \cdot 3}^{\ 3} E_{3 \cdot 1}^{\ 1} - E_{3 \cdot 3}^{\ 1} E_{3 \cdot 1}^{\ 3} = (-1)0 - 0 = 0$$
$$E_{3 \cdot 3}^{\ 2} = E_{3 \cdot 1}^{\ 3} E_{3 \cdot 2}^{\ 1} - E_{3 \cdot 1}^{\ 1} E_{3 \cdot 2}^{\ 3} = 0 * 1 - 0 = 0$$
$$E_{3 \cdot 1}^{\ 3} = E_{3 \cdot 2}^{\ 1} E_{3 \cdot 3}^{\ 2} - E_{3 \cdot 2}^{\ 2} E_{3 \cdot 3}^{\ 1} = 1 * 0 - 0 = 0$$
$$E_{3 \cdot 2}^{\ 3} = E_{3 \cdot 3}^{\ 1} E_{3 \cdot 1}^{\ 2} - E_{3 \cdot 3}^{\ 2} E_{3 \cdot 1}^{\ 1} = 0 * 1 - 0 * 1 = 0$$
$$E_{3 \cdot 3}^{\ 3} = E_{3 \cdot 1}^{\ 1} E_{3 \cdot 2}^{\ 2} - E_{3 \cdot 1}^{\ 2} E_{3 \cdot 2}^{\ 1} = 1 * 0 - 1 * 1 = -1$$

□

*Remark* 8.2. It can be verified directly that $\overline{E}_1$ is linear automorphism. Let
$$a = a^0 + a^1 i + a^2 j + a^3 k$$
$$b = b^0 + b^1 i + b^2 j + b^3 k$$
Then
$$ab = a^0 b^0 - a^1 b^1 - a^2 b^2 - a^3 b^3 + (a^0 b^1 + a^1 b^0 + a^2 b^3 - a^3 b^2)i$$
$$+ (a^0 b^2 + a^2 b^0 + a^3 b^1 - a^1 b^3)j + (a^0 b^3 + a^3 b^0 + a^1 b^2 - a^2 b^1)k$$
$$\overline{E}_1 \circ a = a^0 + a^2 i + a^3 j + a^1 k$$
$$\overline{E}_1 \circ b = b^0 + b^2 i + b^3 j + b^1 k$$
$$(\overline{E}_1 \circ a)(\overline{E}_1 \circ b) = a^0 b^0 - a^2 b^2 - a^3 b^3 - a^1 b^1 + (a^0 b^2 + a^2 b^0 + a^3 b^1 - a^1 b^3)i$$
$$+ (a^0 b^3 + a^3 b^0 + a^1 b^2 - a^2 b^1)j + (a^0 b^1 + a^1 b^0 + a^2 b^3 - a^3 b^2)k$$
$$= \overline{E}_1 \circ (ab)$$
□

**Theorem 8.3.** *The mapping $\overline{E}_3$ has form*
$$\overline{E}_3 \circ a = -\frac{1}{2}(iai + jaj + iaj + jai)$$

*Proof.* According to the theorem [3]-3.3.4, standard components of the mapping $\overline{E}_3$ have form
$$E_1^{11} = -\tfrac{1}{2} \quad E_1^{22} = -\tfrac{1}{2} \quad E_1^{21} = -\tfrac{1}{2} \quad E_1^{12} = -\tfrac{1}{2}$$
□

**Theorem 8.4.** *We can identify the mapping*
$$a \circ \overline{E}_3 : H \to H \quad a \in H$$
*and matrix*

(8.1)
$$J_{3l \cdot a} = \begin{pmatrix} a^0 & -a^2 & -a^1 & a^3 \\ a^1 & -a^3 & a^0 & -a^2 \\ a^2 & a^0 & a^3 & a^1 \\ a^3 & a^1 & -a^2 & -a^0 \end{pmatrix}$$

(8.2)
$$J_{3l \cdot a} = J_{l \cdot a} \circ E_3$$

*Proof.* The product of quaternions
$$a = a^0 + a^1 i + a^2 j + a^3 k$$
and
$$\overline{E}_3 \circ x = x^0 + x^2 i + x^1 j - x^3 k$$
has form
$$a \circ \overline{E}_3 \circ x = a^0 x^0 - a^1 x^2 - a^2 x^1 + a^3 x^3 + (a^0 x^2 + a^1 x^0 - a^2 x^3 - a^3 x^1)i$$
$$+ (a^0 x^1 + a^2 x^0 + a^3 x^2 + a^1 x^3)j + (-a^0 x^3 + a^3 x^0 + a^1 x^1 - a^2 x^2)k$$

Therefore, function $f_a \circ x = a \circ \overline{E}_1 \circ x$ has Jacobian matrix (8.1). The equation (8.2) follows from the chain of equations

$$J_{l \cdot a} \circ E_3 = \begin{pmatrix} a^0 & -a^1 & -a^2 & -a^3 \\ a^1 & a^0 & -a^3 & a^2 \\ a^2 & a^3 & a^0 & -a^1 \\ a^3 & -a^2 & a^1 & a^0 \end{pmatrix} \circ \begin{pmatrix} 1 & 0 & 0 & 0 \\ 0 & 0 & 1 & 0 \\ 0 & 1 & 0 & 0 \\ 0 & 0 & 0 & -1 \end{pmatrix}$$

$$= \begin{pmatrix} a^0 & -a^2 & -a^1 & a^3 \\ a^1 & -a^3 & a^0 & -a^2 \\ a^2 & a^0 & a^3 & a^1 \\ a^3 & a^1 & -a^2 & -a^0 \end{pmatrix} = J_{3l \cdot a}$$

□

## 9. Mapping $\overline{I}$

Antilinear automorphism

$$\overline{I} : H \to H \quad \overline{I} \circ x = x^*$$

$$I_0^0 = 1 \quad I_1^1 = -1 \quad I_2^2 = -1 \quad I_3^3 = -1$$

$$I = \begin{pmatrix} 1 & 0 & 0 & 0 \\ 0 & -1 & 0 & 0 \\ 0 & 0 & -1 & 0 \\ 0 & 0 & 0 & -1 \end{pmatrix}$$

is called conjugation of quaternion algebra.

**Theorem 9.1.** *We can identify the mapping*

$$a \circ \overline{I} : H \to H \quad a \in H$$

*and matrix*

(9.1) $$I_{l \cdot a} = \begin{pmatrix} a^0 & a^1 & a^2 & a^3 \\ a^1 & -a^0 & a^3 & -a^2 \\ a^2 & -a^3 & -a^0 & a^1 \\ a^3 & a^2 & -a^1 & -a^0 \end{pmatrix}$$

(9.2) $$I_{l \cdot a} = J_{l \cdot a} \circ I$$

*Proof.* The product of quaternions

$$a = a^0 + a^1 i + a^2 j + a^3 k$$

and

$$I \circ x = x^0 - x^1 i - x^2 j - x^3 k$$

has form

$$a \circ I \circ x = a^0 x^0 + a^1 x^1 + a^2 x^2 + a^3 x^3 + (-a^0 x^1 + a^1 x^0 - a^2 x^3 + a^3 x^2)i$$
$$+ (-a^0 x^2 + a^2 x^0 - a^3 x^1 + a^1 x^3)j + (-a^0 x^3 + a^3 x^0 - a^1 x^2 + a^2 x^1)k$$

Therefore, function $f_a \circ x = a \circ I \circ x$ has Jacobian matrix (9.1). The equation (9.2) follows from the chain of equations

$$J_{l \cdot a} \circ I = \begin{pmatrix} a^0 & -a^1 & -a^2 & -a^3 \\ a^1 & a^0 & -a^3 & a^2 \\ a^2 & a^3 & a^0 & -a^1 \\ a^3 & -a^2 & a^1 & a^0 \end{pmatrix} \circ \begin{pmatrix} 1 & 0 & 0 & 0 \\ 0 & -1 & 0 & 0 \\ 0 & 0 & -1 & 0 \\ 0 & 0 & 0 & -1 \end{pmatrix}$$

$$= \begin{pmatrix} a^0 & a^1 & a^2 & a^3 \\ a^1 & -a^0 & a^3 & -a^2 \\ a^2 & -a^3 & -a^0 & a^1 \\ a^3 & a^2 & -a^1 & -a^0 \end{pmatrix} = I_{l \cdot a}$$

□

**Theorem 9.2.** *We can identify the mapping*

$$a \star \overline{I} : H \to H \quad a \in H$$

*and matrix*

(9.3) $$I_{r \cdot a} = \begin{pmatrix} a^0 & a^1 & a^2 & a^3 \\ a^1 & -a^0 & -a^3 & a^2 \\ a^2 & a^3 & -a^0 & -a^1 \\ a^3 & -a^2 & a^1 & -a^0 \end{pmatrix}$$

(9.4) $$I_{r \cdot a} = J_{r \cdot a} \circ I$$

*Proof.* The product of quaternions

$$I \circ x = x^0 - x^1 i - x^2 j - x^3 k$$

and

$$a = a^0 + a^1 i + a^2 j + a^3 k$$

has form

$$(I \circ x)a = x^0 a^0 + x^1 a^1 + x^2 a^2 + x^3 a^3 + (x^0 a^1 - x^1 a^0 - x^2 a^3 + x^3 a^2)i$$
$$+ (x^0 a^2 - x^2 a^0 - x^3 a^1 + x^1 a^3)j + (x^0 a^3 - x^3 a^0 - x^1 a^2 + x^2 a^1)k$$

Therefore, function $f_a \circ x = a \circ I \circ x$ has Jacobian matrix (9.3). The equation (9.4) follows from the chain of equations

$$J_{r \cdot a} \circ I = \begin{pmatrix} a^0 & -a^1 & -a^2 & -a^3 \\ a^1 & a^0 & a^3 & -a^2 \\ a^2 & -a^3 & a^0 & a^1 \\ a^3 & a^2 & -a^1 & a^0 \end{pmatrix} \circ \begin{pmatrix} 1 & 0 & 0 & 0 \\ 0 & -1 & 0 & 0 \\ 0 & 0 & -1 & 0 \\ 0 & 0 & 0 & -1 \end{pmatrix}$$

$$= \begin{pmatrix} a^0 & a^1 & a^2 & a^3 \\ a^1 & -a^0 & -a^3 & a^2 \\ a^2 & a^3 & -a^0 & -a^1 \\ a^3 & -a^2 & a^1 & -a^0 \end{pmatrix} = I_{r \cdot a}$$

□

## 10. Mapping $\overline{I}_1$

**Theorem 10.1.** *Since the mapping*

$$f_2 : H \to H$$

*is antilinear automorphism of quaternion algebra and the mapping*

$$f_1 : H \to H$$

*is linear automorphism of quaternion algebra, then the mapping $f_2 \circ f_1$ is antilinear automorphism of quaternion algebra.*

*Proof.* Since quaternion algebra is $R$-algebra, $R \subset \operatorname{Re} H$, then, for $r \in R$, $r^* = r$. Therefore,

$$f_2 \circ f_1 \circ (ra) = f_2 \circ (rf_1 \circ a) = (f_2 \circ f_1 \circ a)r^*$$

The statement of the theorem is corollary of following equations

$$f_2 \circ f_1 \circ (a+b) = f_2 \circ (f_1 \circ a + f_1 \circ b) = f_2 \circ f_1 \circ a + f_2 \circ f_1 \circ b$$
$$f_2 \circ f_1 \circ (ab) = f_2((f_1 \circ a)(f_1 \circ b)) = (f_2 \circ f_1 \circ b)(f_2 \circ f_1 \circ a)$$

□

**Theorem 10.2.** *The mapping*

$$\overline{I}_1 : H \to H \quad \overline{I}_1 \circ x = x^0 - x^2 \overline{e}_1 - x^3 \overline{e}_2 - x^1 \overline{e}_3$$

$$I_{1 \cdot 0}^0 = 1 \quad I_{1 \cdot 2}^1 = -1 \quad I_{1 \cdot 3}^2 = -1 \quad I_{1 \cdot 1}^3 = -1$$

$$I_1 = \begin{pmatrix} 1 & 0 & 0 & 0 \\ 0 & 0 & -1 & 0 \\ 0 & 0 & 0 & -1 \\ 0 & -1 & 0 & 0 \end{pmatrix}$$

*is the antilinear automorphism of quaternion algebra.*

*Proof.* The statement of the theorem is corollary of the theorem 10.1 and the equation
$$\overline{I}_1 = \overline{I} \circ \overline{E}_1$$

□

**Theorem 10.3.** *We can identify the mapping*
$$a \circ \overline{I}_1 : H \to H \quad a \in H$$
*and matrix*

(10.1) $$I_{1l \cdot a} = \begin{pmatrix} a^0 & a^3 & a^1 & a^2 \\ a^1 & -a^2 & -a^0 & a^3 \\ a^2 & a^1 & -a^3 & -a^0 \\ a^3 & -a^0 & a^2 & -a^1 \end{pmatrix}$$

(10.2) $$I_{1l \cdot a} = J_{1l \cdot a} \circ I_1$$

*Proof.* The product of quaternions
$$a = a^0 + a^1 i + a^2 j + a^3 k$$
and
$$\overline{I}_1 \circ x = x^0 - x^2 i - x^3 j - x^1 k$$
has form
$$a \circ \overline{I}_1 \circ x = a^0 x^0 + a^1 x^2 + a^2 x^3 + a^3 x^1 + (-a^0 x^2 + a^1 x^0 - a^2 x^1 + a^3 x^3)i$$
$$+ (-a^0 x^3 + a^2 x^0 - a^3 x^2 + a^1 x^1)j + (-a^0 x^1 + a^3 x^0 - a^1 x^3 + a^2 x^2)k$$
Therefore, function $f_a \circ x = a \circ I \circ x$ has Jacobian matrix (10.1). The equation (10.2) follows from the chain of equations

$$J_{l \cdot a} \circ I_1 = \begin{pmatrix} a^0 & -a^1 & -a^2 & -a^3 \\ a^1 & a^0 & -a^3 & a^2 \\ a^2 & a^3 & a^0 & -a^1 \\ a^3 & -a^2 & a^1 & a^0 \end{pmatrix} \circ \begin{pmatrix} 1 & 0 & 0 & 0 \\ 0 & 0 & -1 & 0 \\ 0 & 0 & 0 & -1 \\ 0 & -1 & 0 & 0 \end{pmatrix}$$

$$= \begin{pmatrix} a^0 & a^3 & a^1 & a^2 \\ a^1 & -a^2 & -a^0 & a^3 \\ a^2 & a^1 & -a^3 & -a^0 \\ a^3 & -a^0 & a^2 & -a^1 \end{pmatrix} = I_{1l \cdot a}$$

□

**Theorem 10.4.** *We can identify the mapping*
$$a \star \overline{I}_1 : H \to H \quad a \in H$$

*and matrix*

$$(10.3) \qquad I_{1r\cdot a} = \begin{pmatrix} a^0 & a^3 & a^1 & a^2 \\ a^1 & a^2 & -a^0 & -a^3 \\ a^2 & -a^1 & a^3 & -a^0 \\ a^3 & -a^0 & -a^2 & a^1 \end{pmatrix}$$

$$(10.4) \qquad I_{1r\cdot a} = J_{1r\cdot a}\circ I_1$$

*Proof.* The product of quaternions

$$\overline{I}_1 \circ x = x^0 - x^2 i - x^3 j - x^1 k$$

and

$$a = a^0 + a^1 i + a^2 j + a^3 k$$

has form

$$(I \circ x)a = x^0 a^0 + x^2 a^1 + x^3 a^2 + x^1 a^3 + (x^0 a^1 - x^2 a^0 - x^3 a^3 + x^1 a^2)i$$
$$+ (x^0 a^2 - x^3 a^0 - x^1 a^1 + x^2 a^3)j + (x^0 a^3 - x^1 a^0 - x^2 a^2 + x^3 a^1)k$$

Therefore, function $f_a \circ x = a \star I \circ x$ has Jacobian matrix (10.3). The equation (10.4) follows from the chain of equations

$$J_{r\cdot a}\circ I_1 = \begin{pmatrix} a^0 & -a^1 & -a^2 & -a^3 \\ a^1 & a^0 & a^3 & -a^2 \\ a^2 & -a^3 & a^0 & a^1 \\ a^3 & a^2 & -a^1 & a^0 \end{pmatrix} \circ \begin{pmatrix} 1 & 0 & 0 & 0 \\ 0 & 0 & -1 & 0 \\ 0 & 0 & 0 & -1 \\ 0 & -1 & 0 & 0 \end{pmatrix}$$

$$= \begin{pmatrix} a^0 & a^3 & a^1 & a^2 \\ a^1 & a^2 & -a^0 & -a^3 \\ a^2 & -a^1 & a^3 & -a^0 \\ a^3 & -a^0 & -a^2 & a^1 \end{pmatrix} = I_{1r\cdot a}$$

□

## 11. Mapping $\overline{I}_2$

**Theorem 11.1.** *The mapping*

$$\overline{I}_2 : H \to H \qquad \overline{I}_2 \circ x = x^0 - x^3 \overline{e}_1 - x^1 \overline{e}_2 - x^2 \overline{e}_3$$

$$I_{2\cdot 0}^0 = 1 \quad I_{2\cdot 3}^1 = -1 \quad I_{2\cdot 1}^2 = -1 \quad I_{2\cdot 2}^3 = -1$$

$$I_2 = \begin{pmatrix} 1 & 0 & 0 & 0 \\ 0 & 0 & 0 & -1 \\ 0 & -1 & 0 & 0 \\ 0 & 0 & -1 & 0 \end{pmatrix}$$

*is the antilinear automorphism of quaternion algebra.*

*Proof.* The statement of the theorem is corollary of the theorem 10.1 and the equation
$$\overline{I}_2 = \overline{I} \circ \overline{E}_2$$
□

**Theorem 11.2.** *We can identify the mapping*
$$a \circ \overline{I}_2 : H \to H \quad a \in H$$
*and matrix*

(11.1)
$$I_{2l \cdot a} = \begin{pmatrix} a^0 & a^2 & a^3 & a^1 \\ a^1 & a^3 & -a^2 & -a^0 \\ a^2 & -a^0 & a^1 & -a^3 \\ a^3 & -a^1 & -a^0 & a^2 \end{pmatrix}$$

(11.2)
$$I_{2l \cdot a} = J_{l \cdot a} \circ I_2$$

*Proof.* The product of quaternions
$$a = a^0 + a^1 i + a^2 j + a^3 k$$
and
$$\overline{I}_2 \circ x = x^0 - x^3 i - x^1 j - x^2 k$$
has form
$$a \circ \overline{I}_2 \circ x = a^0 x^0 + a^1 x^3 + a^2 x^1 + a^3 x^2 + (-a^0 x^3 + a^1 x^0 - a^2 x^2 + a^3 x^1)i$$
$$+ (-a^0 x^1 + a^2 x^0 - a^3 x^3 + a^1 x^2)j + (-a^0 x^2 + a^3 x^0 - a^1 x^1 + a^2 x^3)k$$

Therefore, function $f_a \circ x = a \circ \overline{I}_2 \circ x$ has Jacobian matrix (11.1). The equation (11.2) follows from the chain of equations

$$J_{l \cdot a} \circ I_2 = \begin{pmatrix} a^0 & -a^1 & -a^2 & -a^3 \\ a^1 & a^0 & -a^3 & a^2 \\ a^2 & a^3 & a^0 & -a^1 \\ a^3 & -a^2 & a^1 & a^0 \end{pmatrix} \circ \begin{pmatrix} 1 & 0 & 0 & 0 \\ 0 & 0 & 0 & -1 \\ 0 & -1 & 0 & 0 \\ 0 & 0 & -1 & 0 \end{pmatrix}$$

$$= \begin{pmatrix} a^0 & a^2 & a^3 & a^1 \\ a^1 & a^3 & -a^2 & -a^0 \\ a^2 & -a^0 & a^1 & -a^3 \\ a^3 & -a^1 & -a^0 & a^2 \end{pmatrix} = I_{2l \cdot a}$$

□

## 12. Linear Mapping

**Theorem 12.1.** *Linear mapping of quaternion algebra*

$$\overline{f} : H \to H$$

*has unique expansion*

(12.1) $$\overline{f} = a_0 \circ \overline{E} + a_1 \circ \overline{E}_1 + a_2 \circ \overline{E}_3 + a_3 \circ \overline{I}$$

*In this case*

(12.2) $$a_0^0 = \frac{1}{2}(f_0^0 + f_1^1 - f_2^1 + f_3^0)$$

(12.3) $$a_0^1 = \frac{1}{2}(f_0^1 - f_1^0 + f_2^0 + f_3^1)$$

(12.4) $$a_0^2 = \frac{1}{2}(f_0^2 - f_1^3 + f_2^3 + f_3^2)$$

(12.5) $$a_0^3 = \frac{1}{2}(f_0^3 + f_1^2 - f_2^2 + f_3^3)$$

(12.6) $$a_1^0 = \frac{1}{2}(-f_1^0 - f_1^1 + f_2^0 + f_2^1 - f_3^0 + f_3^1 + f_3^2 + f_3^3)$$

(12.7) $$a_1^1 = \frac{1}{2}(f_1^0 - f_1^1 - f_2^0 + f_2^1 - f_3^0 - f_3^1 - f_3^2 + f_3^3)$$

(12.8) $$a_1^2 = \frac{1}{2}(-f_1^2 + f_1^3 + f_2^2 - f_2^3 - f_3^0 + f_3^1 - f_3^2 - f_3^3)$$

(12.9) $$a_1^3 = \frac{1}{2}(-f_1^2 - f_1^3 + f_2^2 + f_2^3 - f_3^0 - f_3^1 + f_3^2 - f_3^3)$$

(12.10) $$a_2^0 = \frac{1}{2}(f_1^0 + f_1^2 - f_2^0 + f_2^1 - f_3^1 - f_3^2)$$

(12.11) $$a_2^1 = \frac{1}{2}(f_1^1 + f_1^3 - f_2^0 - f_2^1 + f_3^0 - f_3^3)$$

(12.12) $$a_2^2 = \frac{1}{2}(-f_1^0 + f_1^2 - f_2^2 - f_2^3 + f_3^0 + f_3^3)$$

(12.13) $$a_2^3 = \frac{1}{2}(-f_1^1 + f_1^3 + f_2^2 - f_2^3 + f_3^1 - f_3^2)$$

(12.14) $$a_3^0 = \frac{1}{2}(f_0^0 - f_1^2 - f_2^1 - f_3^3)$$

(12.15) $$a_3^1 = \frac{1}{2}(f_0^1 - f_1^3 + f_2^0 + f_3^2)$$

(12.16) $$a_3^2 = \frac{1}{2}(f_0^2 + f_1^0 + f_2^3 - f_3^1)$$

(12.17) $$a_3^3 = \frac{1}{2}(f_0^3 + f_1^1 - f_2^2 + f_3^0)$$

*Proof.* Linear mapping (12.1) has matrix

$$
(12.18) \quad
\begin{pmatrix}
a_0^0 & -a_0^1 & -a_0^2 & -a_0^3 \\
a_0^1 & a_0^0 & -a_0^3 & a_0^2 \\
a_0^2 & a_0^3 & a_0^0 & -a_0^1 \\
a_0^3 & -a_0^2 & a_0^1 & a_0^0
\end{pmatrix}
+
\begin{pmatrix}
a_1^0 & -a_1^3 & -a_1^1 & -a_1^2 \\
a_1^1 & a_1^2 & a_1^0 & -a_1^3 \\
a_1^2 & -a_1^1 & a_1^3 & a_1^0 \\
a_1^3 & a_1^0 & -a_1^2 & a_1^1
\end{pmatrix}
$$
$$
+
\begin{pmatrix}
a_2^0 & -a_2^2 & -a_2^1 & a_2^3 \\
a_2^1 & -a_2^3 & a_2^0 & -a_2^2 \\
a_2^2 & a_2^0 & a_2^3 & a_2^1 \\
a_2^3 & a_2^1 & -a_2^2 & -a_2^0
\end{pmatrix}
+
\begin{pmatrix}
a_3^0 & a_3^1 & a_3^2 & a_3^3 \\
a_3^1 & -a_3^0 & a_3^3 & -a_3^2 \\
a_3^2 & -a_3^3 & -a_3^0 & a_3^1 \\
a_3^3 & a_3^2 & -a_3^1 & -a_3^0
\end{pmatrix}
$$

From a comparison of matrix of the mapping $\overline{f}$ and matrix (12.18), we get the system of linear equations

(12.19) $\quad f_0^0 = a_0^0 + a_1^0 + a_2^0 + a_3^0$ $\qquad$ (12.27) $\quad f_0^2 = a_0^2 + a_1^2 + a_2^2 + a_3^2$

(12.20) $\quad f_1^0 = -a_0^1 - a_1^3 - a_2^2 + a_3^1$ $\qquad$ (12.28) $\quad f_1^2 = a_0^3 - a_1^1 + a_2^0 - a_3^3$

(12.21) $\quad f_2^0 = -a_0^2 - a_1^1 - a_2^1 + a_3^2$ $\qquad$ (12.29) $\quad f_2^2 = a_0^0 + a_1^3 + a_2^3 - a_3^0$

(12.22) $\quad f_3^0 = -a_0^3 - a_1^2 + a_2^3 + a_3^3$ $\qquad$ (12.30) $\quad f_3^2 = -a_0^1 + a_1^0 + a_2^1 + a_3^1$

(12.23) $\quad f_0^1 = a_0^1 + a_1^1 + a_2^1 + a_3^1$ $\qquad$ (12.31) $\quad f_0^3 = a_0^3 + a_1^3 + a_2^3 + a_3^3$

(12.24) $\quad f_1^1 = a_0^0 + a_1^2 - a_2^3 - a_3^0$ $\qquad$ (12.32) $\quad f_1^3 = -a_0^2 + a_1^0 + a_2^1 + a_3^2$

(12.25) $\quad f_2^1 = -a_0^3 + a_1^0 + a_2^0 + a_3^3$ $\qquad$ (12.33) $\quad f_2^3 = a_0^1 - a_1^2 - a_2^2 - a_3^1$

(12.26) $\quad f_3^1 = a_0^2 - a_1^3 - a_2^2 - a_3^2$ $\qquad$ (12.34) $\quad f_3^3 = a_0^0 + a_1^1 - a_2^0 - a_3^0$

To solve this system of linear equations, I have wrote the software using C#. It is easy verify directly this solution. $\square$

It follows from the theorem [3]-2.6.4 that the set of linear endomorphisms $\mathcal{L}(H;H)$ of quaternion algebra $H$ is isomorphic to tensor product $H \otimes H$. The theorem 12.1 states that we can consider the module $\mathcal{L}(H;H)$ as $H\star$-vector space with basis

$$(\overline{E}, \overline{E}_1, \overline{E}_3, \overline{I})$$

**Example 12.2.** According to the theorem 7.4, mapping

$$a \circ \overline{E}_2 : H \to H \quad a \in H$$

has matrix

$$
\begin{pmatrix}
a^0 & -a^2 & -a^3 & -a^1 \\
a^1 & -a^3 & a^2 & a^0 \\
a^2 & a^0 & -a^1 & a^3 \\
a^3 & a^1 & a^0 & -a^2
\end{pmatrix}
$$

According to the theorem 12.1,

$$a_0^0 = \frac{1}{2}(a^0 + (-a^3) - a^2 + (-a^1))$$

$$a_0^1 = \frac{1}{2}(a^1 - (-a^2) + (-a^3) + a^0)$$

$$a_0^2 = \frac{1}{2}(a^2 - a^1 + a^0 + a^3)$$

$$a_0^3 = \frac{1}{2}(a^3 + a^0 - (-a^1) + (-a^2))$$

$$a_1^0 = \frac{1}{2}(-(-a^2) - (-a^3) + (-a^3) + a^2 - (-a^1) + a^0 + a^3 + (-a^2))$$

$$a_1^1 = \frac{1}{2}((-a^2) - (-a^3) - (-a^3) + a^2 - (-a^1) - a^0 - a^3 + (-a^2))$$

$$a_1^2 = \frac{1}{2}(-a^0 + a^1 + (-a^1) - a^0 - (-a^1) + a^0 - a^3 - (-a^2))$$

$$a_1^3 = \frac{1}{2}(-a^0 - a^1 + (-a^1) + a^0 - (-a^1) - a^0 + a^3 - (-a^2))$$

$$a_2^0 = \frac{1}{2}((-a^2) + a^0 - (-a^3) + a^2 - a^0 - a^3)$$

$$a_2^1 = \frac{1}{2}((-a^3) + a^1 - (-a^3) - a^2 + (-a^1) - (-a^2))$$

$$a_2^2 = \frac{1}{2}(-(-a^2) + a^0 - (-a^1) - a^0 + (-a^1) + (-a^2))$$

$$a_2^3 = \frac{1}{2}(-(-a^3) + a^1 + (-a^1) - a^0 + a^0 - a^3)$$

$$a_3^0 = \frac{1}{2}(a^0 - a^0 - a^2 - (-a^2))$$

$$a_3^1 = \frac{1}{2}(a^1 - a^1 + (-a^3) + a^3)$$

$$a_3^2 = \frac{1}{2}(a^2 + (-a^2) + a^0 - a^0)$$

$$a_3^3 = \frac{1}{2}(a^3 + (-a^3) - (-a^1) + (-a^1))$$

Therefore, $a_2 = a_3 = 0$. □

**Example 12.3.** According to the theorem 5.2, mapping

$$a \star \overline{E} : H \to H \quad a \in H$$

has matrix

$$\begin{pmatrix} a^0 & -a^1 & -a^2 & -a^3 \\ a^1 & a^0 & a^3 & -a^2 \\ a^2 & -a^3 & a^0 & a^1 \\ a^3 & a^2 & -a^1 & a^0 \end{pmatrix}$$

According to the theorem 12.1,

$$a_0^0 = \frac{1}{2}(a^0 + a^0 - a^3 + (-a^3))$$

$$a_0^1 = \frac{1}{2}(a^1 - (-a^1) + (-a^2) + (-a^2))$$

$$a_0^2 = \frac{1}{2}(a^2 - a^2 + (-a^1) + a^1)$$

$$a_0^3 = \frac{1}{2}(a^3 + (-a^3) - a^0 + a^0)$$

$$a_1^0 = \frac{1}{2}(-(-a^1) - a^0 + (-a^2) + a^3 - (-a^3) + (-a^2) + a^1 + a^0)$$

$$a_1^1 = \frac{1}{2}((-a^1) - a^0 - (-a^2) + a^3 - (-a^3) - (-a^2) - a^1 + a^0)$$

$$a_1^2 = \frac{1}{2}(-(-a^3) + a^2 + a^0 - (-a^1) - (-a^3) + (-a^2) - a^1 - a^0)$$

$$a_1^3 = \frac{1}{2}(-(-a^3) - a^2 + a^0 + (-a^1) - (-a^3) - (-a^2) + a^1 - a^0)$$

$$a_2^0 = \frac{1}{2}((-a^1) + (-a^3) - (-a^2) + a^3 - (-a^2) - a^1)$$

$$a_2^1 = \frac{1}{2}(a^0 + a^2 - (-a^2) - a^3 + (-a^3) - a^0)$$

$$a_2^2 = \frac{1}{2}(-(-a^1) + (-a^3) - a^0 - (-a^1) + (-a^3) + a^0)$$

$$a_2^3 = \frac{1}{2}(-a^0 + a^2 + a^0 - (-a^1) + (-a^2) - a^1)$$

$$a_3^0 = \frac{1}{2}(a^0 - (-a^3) - a^3 - a^0)$$

$$a_3^1 = \frac{1}{2}(a^1 - a^2 + (-a^2) + a^1)$$

$$a_3^2 = \frac{1}{2}(a^2 + (-a^1) + (-a^1) - (-a^2))$$

$$a_3^3 = \frac{1}{2}(a^3 + a^0 - a^0 + (-a^3))$$

This is very important, $a_3 \neq 0$.　　　　　　　　　　　　　　　　　　　　□

**Example 12.4.** According to the theorems 5.1, 5.2, mapping

$$f : H \to H \quad f = a \circ E + a \star E \quad a \in H$$

$$f \circ x = ax + xa$$

has matrix

$$\begin{pmatrix} a^0 & -a^1 & -a^2 & -a^3 \\ a^1 & a^0 & -a^3 & a^2 \\ a^2 & a^3 & a^0 & -a^1 \\ a^3 & -a^2 & a^1 & a^0 \end{pmatrix} + \begin{pmatrix} a^0 & -a^1 & -a^2 & -a^3 \\ a^1 & a^0 & a^3 & -a^2 \\ a^2 & -a^3 & a^0 & a^1 \\ a^3 & a^2 & -a^1 & a^0 \end{pmatrix}$$

$$= \begin{pmatrix} 2a^0 & -2a^1 & -2a^2 & -2a^3 \\ 2a^1 & 2a^0 & 0 & 0 \\ 2a^2 & 0 & 2a^0 & 0 \\ 2a^3 & 0 & 0 & 2a^0 \end{pmatrix}$$

According to the theorem 12.1,

$$a_0^0 = \frac{1}{2}(2a^0 + 2a^0 - 0 + (-2a^3))$$
$$a_0^1 = \frac{1}{2}(2a^1 - (-2a^1) + (-2a^2) + 0)$$
$$a_0^2 = \frac{1}{2}(2a^2 - 0 + 0 + 0)$$
$$a_0^3 = \frac{1}{2}(2a^3 + 0 - 2a^0 + 2a^0)$$

$$a_1^0 = \frac{1}{2}(-(-2a^1) - 2a^0 + (-2a^2) + 0 - (-2a^3) + 0 + 0 + 2a^0)$$
$$a_1^1 = \frac{1}{2}((-2a^1) - 2a^0 - (-2a^2) + 0 - (-2a^3) - 0 - 0 + 2a^0)$$
$$a_1^2 = \frac{1}{2}(-0 + 0 + 2a^0 - 0 - (-2a^3) + 0 - 0 - 2a^0)$$
$$a_1^3 = \frac{1}{2}(-0 - 0 + 2a^0 + 0 - (-2a^3) - 0 + 0 - 2a^0)$$

$$a_2^0 = \frac{1}{2}((-2a^1) + 0 - (-2a^2) + 0 - 0 - 0)$$
$$a_2^1 = \frac{1}{2}(2a^0 + 0 - (-2a^2) - 0 + (-2a^3) - 2a^0)$$
$$a_2^2 = \frac{1}{2}(-(-2a^1) + 0 - 2a^0 - 0 + (-2a^3) + 2a^0)$$
$$a_2^3 = \frac{1}{2}(-2a^0 + 0 + 2a^0 - 0 + 0 - 0)$$

$$a_3^0 = \frac{1}{2}(2a^0 - 0 - 0 - 2a^0)$$
$$a_3^1 = \frac{1}{2}(2a^1 - 0 + (-2a^2) + 0)$$
$$a_3^2 = \frac{1}{2}(2a^2 + (-2a^1) + 0 - 0)$$
$$a_3^3 = \frac{1}{2}(2a^3 + 2a^0 - 2a^0 + (-2a^3))$$

This is very important, $a_3 \neq 0$. □

# 14. Special Symbols and Notations

# Линейные отображения алгебры кватернионов

Александр Клейн

Аннотация. В статье рассмотрены линейные и антилинейные автоморфизмы алгебры кватернионов. Доказана теорема о единственности разложения $R$-линейного отображения алгебры кватернионов относительно заданного множества линейных и антилинейных автоморфизмов.

Содержание

## 1. Предисловие к изданию 1

По своему характеру эта статья напоминает справочник элементарных функций. Однако центральное место в статье занимает теорема 12.1, которая является мощным источником новых идей.

Когда я писал статью [2], до последней минуты я надеялся, что мне удастся доказать, что $D$-линейное отображение в $D$-алгебре $A$ с сопряжением можно разложить в сумму $A$-линейного и $A$-антилинейного отображений. Так как для решения этой задачи необходимо составить систему линейных уравнений, то я в конце концов понял что даже в случае алгебры кватернионов эта задача имеет отрицательный ответ.

Aleks_Kleyn@MailAPS.org.
http://sites.google.com/site/AleksKleyn/.
http://arxiv.org/a/kleyn_a_1.
http://AleksKleyn.blogspot.com/.

Утверждение, что в алгебре кватернионов существует нетривиальный линейный автоморфизм, оказалось для меня полной неожиданностью.[1] Поэтому мне потребовалось время, чтобы понять насколько это утверждение важно.

Я начал исследование алгебр со счётным базисом. И вдруг я понял, что если я имею 3 линейных автоморфизма (а именно, $\overline{E}, \overline{E}_1, \overline{E}_2$ ) и антилинейный автоморфизм $\overline{I}$, то есть надежда, что для любого $R$-линейного отображения[2] $f$ существует единственное разложение вида

(1.1) $$f = a_0 \circ \overline{E} + a_1 \circ \overline{E}_1 + a_2 \circ \overline{E}_2 + a_3 \circ \overline{I}$$

Решение казалось простым, однако сопутствующая система линейных уравнений вырожденной. Различные попытки изменить третье слагаемое были безуспешными.

Удача улыбнулась мне, когда я спросил себя, является ли отображение $\overline{E}_3$ линейным автоморфизмом.

## 2. Предисловие к изданию 2

Энтони Садбери обратил моё внимание, что любое преобразование вида

(2.1) $$x \to a\,x\,a^{-1} \quad a \in H$$

является линейным автоморфизмом.[3] В частности, отображение $\overline{E}_1$ имеет вид

$$\overline{E}_1 \circ x = (1 + \sqrt{3}(i + j + k))x(1 + \sqrt{3}(i + j + k))^{-1}$$

Можно ли доказать, что каждый линейный автоморфизм имеет вид (2.1)?

Теорема 12.1 имеет ещё одно, очень важное следствие. В разделе [1]-2.4, я рассмотрел $D$-векторное пространство (в данном случае, $D = H$). Из теоремы 12.1 следует, что размерность $H$-векторного пространства $V$ не совпадает с размерностью $H\star$-векторного пространства $V$.

## 3. Соглашения

**Соглашение 3.1.** *Пусть $A$ - свободная конечно мерная алгебра. При разложении элемента алгебры $A$ относительно базиса $\overline{\overline{e}}$ мы пользуемся одной и той же корневой буквой для обозначения этого элемента и его координат. Однако в алгебре не принято использовать векторные обозначения. В выражении $a^2$ не ясно - это компонента разложения элемента $a$ относительно базиса или это операция возведения в степень. Для облегчения чтения текста мы будем индекс элемента алгебры выделять цветом. Например,*

$$a = a^i \overline{e}_i$$

□

---

[1]Иногда очень трудно воспроизвести реальную последовательность событий. Однако, когда я готовил пример [2]-7.2, я вспомнил диаграмму симметрии произведения в алгебре кватернионов, приведённую в [4], с. 7.

[2]Поскольку $R$-линейные отображения являются основным инструментом в моём исследовании математического анализа, я обычно пользуюсь термином линейные отображения. Однако в тех случаях, когда мы рассматриваем несколько алгебр одновременно, необходимо указать над какой алгеброй отображение линейно.

[3]Подобные отображения хорошо известны (смотри, например, раздел [3]-3.2, а также с. [5]-643). Однако я об этом варианте не думал, так как моё внимание было сосредоточено на отображениях вида $y = af(x)$.

**Соглашение 3.2.** *Если свободная конечномерная алгебра имеет единицу, то мы будем отождествлять вектор базиса $\overline{e}_0$ с единицей алгебры.* □

**Соглашение 3.3.** *Для данной D-алгебры A левый сдвиг $a\circ$ определён равенством*
$$a\circ x = ax$$
*и правый сдвиг $a\star$ определён равенством*
$$a\star x = xa$$

□

Без сомнения, у читателя могут быть вопросы, замечания, возражения. Я буду признателен любому отзыву.

## 4. Линейный автоморфизм алгебры кватернионов

**Теорема 4.1.** *Координаты линейного автоморфизма алгебры кватернионов удовлетворяют системе уравнений*

$$(4.1)\quad \begin{cases} r^1_1 = r^2_2 r^3_3 - r^3_2 r^2_3 & r^1_2 = r^2_3 r^3_1 - r^3_3 r^2_1 & r^1_3 = r^2_1 r^3_2 - r^3_1 r^2_2 \\ r^2_1 = r^3_2 r^1_3 - r^1_2 r^3_3 & r^2_2 = r^3_3 r^1_1 - r^3_1 r^1_3 & r^2_3 = r^3_1 r^1_2 - r^1_1 r^3_2 \\ r^3_1 = r^1_2 r^2_3 - r^2_2 r^1_3 & r^3_2 = r^1_3 r^2_1 - r^2_3 r^1_1 & r^3_3 = r^1_1 r^2_2 - r^2_1 r^1_2 \end{cases}$$

*Доказательство.* Согласно теоремам [3]-3.3.1, [2]-6.4, линейный автоморфизм алгебры кватернионов удовлетворяет уравнениям

$$(4.2)\quad \begin{array}{llll} r^l_0 = r^p_0 r^q_0 C^l_{pq} & r^l_1 = r^p_0 r^q_1 C^l_{pq} & r^l_2 = r^p_0 r^q_2 C^l_{pq} & r^l_3 = r^p_0 r^q_3 C^l_{pq} \\ r^l_1 = r^p_1 r^q_0 C^l_{pq} & -r^l_0 = r^p_1 r^q_1 C^l_{pq} & r^l_3 = r^p_1 r^q_2 C^l_{pq} & -r^l_2 = r^p_1 r^q_3 C^l_{pq} \\ r^l_2 = r^p_2 r^q_0 C^l_{pq} & -r^l_3 = r^p_2 r^q_1 C^l_{pq} & -r^l_0 = r^p_2 r^q_2 C^l_{pq} & r^l_0 = r^p_2 r^q_3 C^l_{pq} \\ r^l_3 = r^p_3 r^q_0 C^l_{pq} & r^l_2 = r^p_3 r^q_1 C^l_{pq} & -r^l_1 = r^p_3 r^q_2 C^l_{pq} & -r^l_0 = r^p_3 r^q_3 C^l_{pq} \end{array}$$

Из равенства (4.2) следует

$$(4.3)\quad \begin{array}{llll} r^l_1 = r^p_0 r^q_1 C^l_{pq} = r^p_2 r^q_3 C^l_{pq} & r^p_0 r^q_1 C^l_{pq} = r^p_0 r^q_1 C^l_{qp} & r^p_2 r^q_3 C^l_{pq} = -r^p_2 r^q_3 C^l_{qp} \\ r^l_2 = r^p_0 r^q_2 C^l_{pq} = r^p_3 r^q_1 C^l_{pq} & r^p_0 r^q_2 C^l_{pq} = r^p_0 r^q_2 C^l_{qp} & r^p_3 r^q_1 C^l_{pq} = -r^p_1 r^q_3 C^l_{pq} \\ r^l_3 = r^p_0 r^q_3 C^l_{pq} = r^p_1 r^q_2 C^l_{pq} & r^p_0 r^q_3 C^l_{pq} = r^p_0 r^q_3 C^l_{qp} & r^p_1 r^q_2 C^l_{pq} = -r^p_1 r^q_2 C^l_{qp} \end{array}$$

$$(4.4)\quad r^l_0 = r^p_0 r^q_0 C^l_{pq} = -r^p_1 r^q_1 C^l_{pq} = -r^p_2 r^q_2 C^l_{pq} = -r^p_3 r^q_3 C^l_{pq}$$

Если $l=0$, то из равенства
$$C^0_{pq} = C^0_{qp}$$
следует
$$(4.5)\quad r^p_i r^q_j C^0_{pq} = r^p_i r^q_j C^0_{qp}$$

Из равенства (4.3) для $l=0$ и равенства (4.5), следует
$$(4.6)\quad r^0_1 = r^0_2 = r^0_3 = 0$$

Если $l = 1, 2, 3$, то равенство (4.3) можно записать в виде

$$
(4.7)\begin{cases}
r_i^l = r_0^l r_i^0 C_{l0}^l + r_0^0 r_i^l C_{0l}^l + r_0^a r_i^b C_{ab}^l + r_0^b r_i^a C_{ba}^l \\
\quad r_0^l r_i^0 C_{l0}^l + r_0^0 r_i^l C_{0l}^l + r_0^a r_i^b C_{ab}^l + r_0^b r_i^a C_{ba}^l \\
= r_0^l r_i^0 C_{0l}^l + r_0^0 r_i^l C_{l0}^l + r_0^a r_i^b C_{ba}^l + r_0^b r_i^a C_{ab}^l \\
\quad i = 1, 2, 3 \\
r_i^l = r_k^0 r_j^l C_{0l}^l + r_k^l r_j^0 C_{l0}^l + r_k^a r_j^b C_{ab}^l + r_k^b r_j^a C_{ba}^l \\
\quad r_k^0 r_j^l C_{0l}^l + r_k^l r_j^0 C_{l0}^l + r_k^a r_j^b C_{ab}^l + r_k^b r_j^a C_{ba}^l \\
= -r_k^0 r_j^l C_{l0}^l - r_k^l r_j^0 C_{0l}^l - r_k^a r_j^b C_{ba}^l - r_k^b r_j^a C_{ab}^l \\
\quad i = 1 \quad k = 2 \quad j = 3 \\
\quad i = 2 \quad k = 3 \quad j = 1 \\
\quad i = 3 \quad k = 1 \quad j = 2 \\
\quad 0 < a < b \quad a \neq l \quad b \neq l
\end{cases}
$$

Из равенств (4.7), (4.6) и равенств

$$
(4.8)\quad \begin{aligned} C_{0l}^l = C_{l0}^l &= 1 \\ C_{ab}^l &= -C_{ba}^l \end{aligned}
$$

следует

$$
(4.9)\begin{cases}
r_i^l = r_0^0 r_i^l + r_0^a r_i^b C_{ab}^l - r_0^b r_i^a C_{ab}^l \\
\quad r_0^0 r_i^l + r_0^a r_i^b C_{ab}^l - r_0^b r_i^a C_{ab}^l \\
= r_0^0 r_i^l - r_0^a r_i^b C_{ab}^l + r_0^b r_i^a C_{ab}^l \\
\quad i = 1, 2, 3 \\
r_i^l = r_k^a r_j^b C_{ab}^l - r_k^b r_j^a C_{ab}^l \\
\quad r_k^a r_j^b C_{ab}^l - r_k^b r_j^a C_{ab}^l \\
= r_k^a r_j^b C_{ab}^l - r_k^b r_j^a C_{ab}^l \\
\quad i = 1 \quad k = 2 \quad j = 3 \\
\quad i = 2 \quad k = 3 \quad j = 1 \\
\quad i = 3 \quad k = 1 \quad j = 2 \\
\quad 0 < a < b \quad a \neq l \quad b \neq l
\end{cases}
$$

Из равенств (4.9) следует

$$
(4.10)\quad\begin{cases} r_i^l = r_0^0 r_i^l \\ \qquad r_0^a r_i^b - r_0^b r_i^a = 0 \\ \qquad i = 1, 2, 3 \\ r_i^l = r_k^a r_j^b C_{ab}^l - r_k^b r_j^a C_{ab}^l \\ \qquad i = 1 \quad k = 2 \quad j = 3 \\ \qquad i = 2 \quad k = 3 \quad j = 1 \\ \qquad i = 3 \quad k = 1 \quad j = 2 \\ \qquad 0 < a < b \quad a \neq l \quad b \neq l \end{cases}
$$

Из равенства (4.10) следует

$$(4.11)\qquad r_0^0 = 1$$

Из равенства (4.4) для $l = 0$ следует

$$
(4.12)\quad\begin{aligned} r_0^0 &= r_0^0 r_0^0 - r_0^1 r_0^1 - r_0^2 r_0^2 - r_0^3 r_0^3 \\ &= -r_i^0 r_i^0 + r_i^1 r_i^1 + r_i^2 r_i^2 + r_i^3 r_i^3 \\ i &= 1, 2, 3 \end{aligned}
$$

Из равенств (4.6), (4.10), (4.12), следует

$$
(4.13)\quad\begin{aligned} 0 &= r_0^1 r_0^1 + r_0^2 r_0^2 + r_0^3 r_0^3 \\ 1 &= r_i^1 r_i^1 + r_i^2 r_i^2 + r_i^3 r_i^3 \\ i &= 1, 2, 3 \end{aligned}
$$

Из равенств (4.13) следует[4]

$$(4.14)\qquad r_0^1 = r_0^2 = r_0^3 = 0$$

Из равенства (4.4) для $l > 0$ следует

$$
(4.15)\quad\begin{aligned} r_0^l &= r_0^l r_0^0 C_{l0}^l + r_0^0 r_0^l C_{0l}^l + r_0^a r_0^b C_{ab}^l + r_0^b r_0^a C_{ba}^l \\ &= -r_i^l r_i^0 C_{l0}^l - r_i^0 r_i^l C_{0l}^l - r_i^a r_i^b C_{ab}^l - r_i^b r_i^a C_{ba}^l \\ i &> 0 \\ l &> 0 \quad 0 < a < b \quad a \neq l \quad b \neq l \end{aligned}
$$

---

[4] Мы здесь опираемся на то, что алгебра кватернионов определена над полем действительных чисел. Если рассматривать алгебру кватернионов над полем комплексных чисел, то уравнение (4.13) определяет конус в комплексном пространстве. Соответственно, у нас шире выбор координат линейного автоморфизма.

Равенства (4.15) тождественно верны в силу равенств (4.6), (4.14), (4.8). Из равенств (4.14), (4.10), следует

$$(4.16) \quad \begin{cases} r_i^l = r_k^a r_j^b C_{ab}^l - r_k^b r_j^a C_{ab}^l \\ \quad i = 1 \quad k = 2 \quad j = 3 \\ \quad i = 2 \quad k = 3 \quad j = 1 \\ \quad i = 3 \quad k = 1 \quad j = 2 \\ \quad l > 0 \quad 0 < a < b \quad a \neq l \quad b \neq l \end{cases}$$

Равенства (4.1) следуют из равенств (4.16). $\square$

## 5. Отображение $\overline{E}$

Очевидно, координаты отображения

$$\overline{E} : H \to H \quad \overline{E} \circ x = x$$

$$E_j^i = \delta_j^i$$

удовлетворяют уравнению (4.1).

**Теорема 5.1.** *Мы можем отождествить отображение*

$$a \circ \overline{E} : H \to H \quad a \in H$$

*и матрицу*

$$(5.1) \quad J_{l \cdot a} = \begin{pmatrix} a^0 & -a^1 & -a^2 & -a^3 \\ a^1 & a^0 & -a^3 & a^2 \\ a^2 & a^3 & a^0 & -a^1 \\ a^3 & -a^2 & a^1 & a^0 \end{pmatrix}$$

*Доказательство.* Произведение кватернионов

$$a = a^0 + a^1 i + a^2 j + a^3 k$$

и

$$x = x^0 + x^1 i + x^2 j + x^3 k$$

имеет вид

$$ax = a^0 x^0 - a^1 x^1 - a^2 x^2 - a^3 x^3 + (a^0 x^1 + a^1 x^0 + a^2 x^3 - a^3 x^2)i$$
$$+ (a^0 x^2 + a^2 x^0 + a^3 x^1 - a^1 x^3)j + (a^0 x^3 + a^3 x^0 + a^1 x^2 - a^2 x^1)k$$

Следовательно, отображение $f_a(x) = ax$ имеет матрицу Якоби (5.1). $\square$

**Теорема 5.2.** *Мы можем отождествить отображение*

$$a \star \overline{E} : H \to H \quad a \in H$$

и матрицу

$$(5.2) \quad J_{r \cdot a} = \begin{pmatrix} a^0 & -a^1 & -a^2 & -a^3 \\ a^1 & a^0 & a^3 & -a^2 \\ a^2 & -a^3 & a^0 & a^1 \\ a^3 & a^2 & -a^1 & a^0 \end{pmatrix}$$

*Доказательство.* Произведение кватернионов

$$x = x^0 + x^1 i + x^2 j + x^3 k$$

и

$$a = a^0 + a^1 i + a^2 j + a^3 k$$

имеет вид

$$xa = x^0 a^0 - x^1 a^1 - x^2 a^2 - x^3 a^3 + (x^0 a^1 + x^1 a^0 + x^2 a^3 - x^3 a^2) i$$
$$+ (x^0 a^2 + x^2 a^0 + x^3 a^1 - x^1 a^3) j + (x^0 a^3 + x^3 a^0 + x^1 a^2 - x^2 a^1) k$$

Следовательно, отображение $f_a(x) = ax$ имеет матрицу Якоби (5.2). □

## 6. Отображение $\overline{E}_1$

**Теорема 6.1.** *Координаты отображения*

$$\overline{E}_1 : H \to H \quad \overline{E}_1 \circ x = x^0 + x^2 i + x^3 j + x^1 k$$

$$E_{1 \cdot 0}^{0} = 1 \quad E_{1 \cdot 2}^{1} = 1 \quad E_{1 \cdot 3}^{2} = 1 \quad E_{1 \cdot 1}^{3} = 1$$

$$E_1 = \begin{pmatrix} 1 & 0 & 0 & 0 \\ 0 & 0 & 1 & 0 \\ 0 & 0 & 0 & 1 \\ 0 & 1 & 0 & 0 \end{pmatrix}$$

*удовлетворяют равенству* (4.1).

*Доказательство.* Мы можем проверить выполнение равенств (4.1) непосредственной проверкой. Однако нетрудно убедиться, что применение перестановки

$$\begin{pmatrix} 1 & 2 & 3 \\ 2 & 3 & 1 \end{pmatrix}$$

на множестве нижних индексов сохраняет множество уравнений. □

*Замечание 6.2.* Мы можем убедиться непосредственной проверкой, что $\overline{E}_1$ - линейный автоморфизм. Пусть

$$a = a^0 + a^1 i + a^2 j + a^3 k$$
$$b = b^0 + b^1 i + b^2 j + b^3 k$$

Тогда

$$ab = a^0 b^0 - a^1 b^1 - a^2 b^2 - a^3 b^3 + (a^0 b^1 + a^1 b^0 + a^2 b^3 - a^3 b^2) i$$
$$+ (a^0 b^2 + a^2 b^0 + a^3 b^1 - a^1 b^3) j + (a^0 b^3 + a^3 b^0 + a^1 b^2 - a^2 b^1) k$$

$$\overline{E}_1 \circ a = a^0 + a^2 i + a^3 j + a^1 k$$
$$\overline{E}_1 \circ b = b^0 + b^2 i + b^3 j + b^1 k$$
$$(\overline{E}_1 \circ a)(\overline{E}_1 \circ b) = a^0 b^0 - a^2 b^2 - a^3 b^3 - a^1 b^1 + (a^0 b^2 + a^2 b^0 + a^3 b^1 - a^1 b^3)i$$
$$+ (a^0 b^3 + a^3 b^0 + a^1 b^2 - a^2 b^1)j + (a^0 b^1 + a^1 b^0 + a^2 b^3 - a^3 b^2)k$$
$$= \overline{E}_1 \circ (ab)$$

$\square$

**Теорема 6.3.** *Отображение* $\overline{E}_1$ *имеет вид*
$$\overline{E}_1 \circ a = \frac{1}{4}(a - iai - jaj - kak - ia + ai - kaj - jak$$
$$-ja - kai + aj - iak - ka - jai - iaj + ak)$$

*Доказательство.* Согласно теореме [3]-3.3.4, стандартные компоненты отображения $\overline{E}_1$ имеют вид

$$\begin{array}{cccc}
E_1^{00} = \frac{1}{4} & E_1^{11} = -\frac{1}{4} & E_1^{22} = -\frac{1}{4} & E_1^{33} = -\frac{1}{4} \\
E_1^{10} = -\frac{1}{4} & E_1^{01} = \frac{1}{4} & E_1^{32} = -\frac{1}{4} & E_1^{23} = -\frac{1}{4} \\
E_1^{20} = -\frac{1}{4} & E_1^{31} = -\frac{1}{4} & E_1^{02} = \frac{1}{4} & E_1^{13} = -\frac{1}{4} \\
E_1^{30} = -\frac{1}{4} & E_1^{21} = -\frac{1}{4} & E_1^{12} = -\frac{1}{4} & E_1^{03} = \frac{1}{4}
\end{array}$$

$\square$

**Теорема 6.4.** *Мы можем отождествить отображение*
$$a \circ \overline{E}_1 : H \to H \quad a \in H$$
*и матрицу*

(6.1) $$J_{1l \cdot a} = \begin{pmatrix} a^0 & -a^3 & -a^1 & -a^2 \\ a^1 & a^2 & a^0 & -a^3 \\ a^2 & -a^1 & a^3 & a^0 \\ a^3 & a^0 & -a^2 & a^1 \end{pmatrix}$$

(6.2) $$J_{1l \cdot a} = J_{l \cdot a} \circ E_1$$

*Доказательство.* Произведение кватернионов
$$a = a^0 + a^1 i + a^2 j + a^3 k$$
и
$$\overline{E}_1 \circ x = x^0 + x^2 i + x^3 j + x^1 k$$
имеет вид
$$a \circ \overline{E}_1 \circ x = a^0 x^0 - a^1 x^2 - a^2 x^3 - a^3 x^1 + (a^0 x^2 + a^1 x^0 + a^2 x^1 - a^3 x^3)i$$
$$+ (a^0 x^3 + a^2 x^0 + a^3 x^2 - a^1 x^1)j + (a^0 x^1 + a^3 x^0 + a^1 x^3 - a^2 x^2)k$$

Следовательно, отображение $f_a \circ x = a \circ \overline{E}_1 \circ x$ имеет матрицу Якоби (6.1). Равенство (6.2) следует из цепочки равенств

$$J_{l \cdot a \circ}{}^\circ E_2 = \begin{pmatrix} a^0 & -a^1 & -a^2 & -a^3 \\ a^1 & a^0 & -a^3 & a^2 \\ a^2 & a^3 & a^0 & -a^1 \\ a^3 & -a^2 & a^1 & a^0 \end{pmatrix} \circ{}^\circ \begin{pmatrix} 1 & 0 & 0 & 0 \\ 0 & 0 & 1 & 0 \\ 0 & 0 & 0 & 1 \\ 0 & 1 & 0 & 0 \end{pmatrix}$$

$$= \begin{pmatrix} a^0 & -a^3 & -a^1 & -a^2 \\ a^1 & a^2 & a^0 & -a^3 \\ a^2 & -a^1 & a^3 & a^0 \\ a^3 & a^0 & -a^2 & a^1 \end{pmatrix} = J_{1l \cdot a}$$

□

## 7. Отображение $\overline{E}_2$

**Теорема 7.1.** *Если отображения*

$$f_1 : H \to H$$
$$f_2 : H \to H$$

*являются линейными автоморфизмами алгебры кватернионов, то отображение $f_2 \circ f_1$ является линейным автоморфизмом алгебры кватернионов.*

*Доказательство.* Утверждение теоремы является следствием следующих равенств

$$f_2 \circ f_1 \circ (a + b) = f_2 \circ (f_1 \circ a + f_1 \circ b) = f_2 \circ f_1 \circ a + f_2 \circ f_1 \circ b$$
$$f_2 \circ f_1 \circ (ra) = f_2 \circ (r(f_1 \circ a)) = r(f_2 \circ f_1 \circ a)$$
$$f_2 \circ f_1 \circ (ab) = f_2((f_1 \circ a)(f_1 \circ b)) = (f_2 \circ f_1 \circ a)(f_2 \circ f_1 \circ b)$$

□

**Теорема 7.2.** *Отображение*

$$\overline{E}_2 : H \to H \quad \overline{E}_2(x) = x^0 + x^3 i + x^1 j + x^2 k$$

$$E_{2 \cdot 0}^0 = 1 \quad E_{2 \cdot 3}^1 = 1 \quad E_{2 \cdot 1}^2 = 1 \quad E_{2 \cdot 2}^3 = 1$$

$$E_2 = \begin{pmatrix} 1 & 0 & 0 & 0 \\ 0 & 0 & 0 & 1 \\ 0 & 1 & 0 & 0 \\ 0 & 0 & 1 & 0 \end{pmatrix}$$

*является линейным автоморфизмом алгебры кватернионов.*

*Доказательство.* Утверждение теоремы является следствием теоремы 7.1 и равенства

$$(E_1)^2 = E_2$$

□

**Теорема 7.3.** *Отображение $\overline{E}_2$ имеет вид*

$$\overline{E}_2 \circ a = \frac{1}{4}(a - iai - jaj - kak + ia - ai - kaj - jak$$
$$+ja - kai - aj - iak + ka - jai - iaj - ak)$$

*Доказательство.* Согласно теореме [3]-3.3.4, стандартные компоненты отображения $\overline{E}_2$ имеют вид

$$\begin{array}{llll}
E_2^{00} = \tfrac{1}{4} & E_2^{11} = -\tfrac{1}{4} & E_2^{22} = -\tfrac{1}{4} & E_2^{33} = -\tfrac{1}{4} \\
E_2^{10} = \tfrac{1}{4} & E_2^{01} = -\tfrac{1}{4} & E_2^{32} = -\tfrac{1}{4} & E_2^{23} = -\tfrac{1}{4} \\
E_2^{20} = \tfrac{1}{4} & E_2^{31} = -\tfrac{1}{4} & E_2^{02} = -\tfrac{1}{4} & E_2^{13} = -\tfrac{1}{4} \\
E_2^{30} = \tfrac{1}{4} & E_2^{21} = -\tfrac{1}{4} & E_2^{12} = -\tfrac{1}{4} & E_2^{03} = -\tfrac{1}{4}
\end{array}$$

$\square$

**Теорема 7.4.** *Мы можем отождествить отображение*

$$a \circ \overline{E}_2 : H \to H \quad a \in H$$

*и матрицу*

(7.1)
$$J_{2l \cdot a} = \begin{pmatrix} a^0 & -a^2 & -a^3 & -a^1 \\ a^1 & -a^3 & a^2 & a^0 \\ a^2 & a^0 & -a^1 & a^3 \\ a^3 & a^1 & a^0 & -a^2 \end{pmatrix}$$

(7.2)
$$J_{2l \cdot a} = J_{l \cdot a} \circ E_2$$

*Доказательство.* Произведение кватернионов

$$a = a^0 + a^1 i + a^2 j + a^3 k$$

и

$$\overline{E}_2 \circ x = x^0 + x^3 i + x^1 j + x^2 k$$

имеет вид

$$a \circ \overline{E}_2 \circ x = a^0 x^0 - a^1 x^3 - a^2 x^1 - a^3 x^2 + (a^0 x^3 + a^1 x^0 + a^2 x^2 - a^3 x^1)i$$
$$+ (a^0 x^1 + a^2 x^0 + a^3 x^3 - a^1 x^2)j + (a^0 x^2 + a^3 x^0 + a^1 x^1 - a^2 x^3)k$$

Следовательно, отображение $f_a \circ x = a \circ \overline{E}_2 \circ x$ имеет матрицу Якоби (7.1). Равенство (7.2) следует из цепочки равенств

$$J_{l \cdot a} \circ E_2 = \begin{pmatrix} a^0 & -a^1 & -a^2 & -a^3 \\ a^1 & a^0 & -a^3 & a^2 \\ a^2 & a^3 & a^0 & -a^1 \\ a^3 & -a^2 & a^1 & a^0 \end{pmatrix} \circ \begin{pmatrix} 1 & 0 & 0 & 0 \\ 0 & 0 & 0 & 1 \\ 0 & 1 & 0 & 0 \\ 0 & 0 & 1 & 0 \end{pmatrix}$$

$$= \begin{pmatrix} a^0 & -a^2 & -a^3 & -a^1 \\ a^1 & -a^3 & a^2 & a^0 \\ a^2 & a^0 & -a^1 & a^3 \\ a^3 & a^1 & a^0 & -a^2 \end{pmatrix} = J_{2l \cdot a}$$

$\square$

## 8. Отображение $\overline{E}_3$

**Теорема 8.1.** *Координаты отображения*

$$\overline{E}_3 : H \to H \quad \overline{E}_3 \circ x = x^0 + x^2 i + x^1 j - x^3 k$$

$$E_{3 \cdot 0}^{0} = 1 \quad E_{3 \cdot 2}^{1} = 1 \quad E_{3 \cdot 1}^{2} = 1 \quad E_{3 \cdot 3}^{3} = -1$$

$$E_3 = \begin{pmatrix} 1 & 0 & 0 & 0 \\ 0 & 0 & 1 & 0 \\ 0 & 1 & 0 & 0 \\ 0 & 0 & 0 & -1 \end{pmatrix}$$

*удовлетворяют равенствам* (4.1).

*Доказательство.* Мы можем проверить выполнение равенств (4.1) непосредственной проверкой

$$E_{3 \cdot 1}^{1} = E_{3 \cdot 2}^{2} E_{3 \cdot 3}^{3} - E_{3 \cdot 2}^{3} E_{3 \cdot 3}^{2} = 0 - 0 = 0$$
$$E_{3 \cdot 2}^{1} = E_{3 \cdot 3}^{2} E_{3 \cdot 1}^{3} - E_{3 \cdot 3}^{3} E_{3 \cdot 1}^{2} = 0 - (-1)1 = 1$$
$$E_{3 \cdot 3}^{1} = E_{3 \cdot 1}^{2} E_{3 \cdot 2}^{3} - E_{3 \cdot 1}^{3} E_{3 \cdot 2}^{2} = 1 * 0 - 0 = 0$$
$$E_{3 \cdot 1}^{2} = E_{3 \cdot 2}^{3} E_{3 \cdot 3}^{1} - E_{3 \cdot 2}^{1} E_{3 \cdot 3}^{3} = 0 - 1(-1) = 1$$
$$E_{3 \cdot 2}^{2} = E_{3 \cdot 3}^{3} E_{3 \cdot 1}^{1} - E_{3 \cdot 3}^{1} E_{3 \cdot 1}^{3} = (-1)0 - 0 = 0$$
$$E_{3 \cdot 3}^{2} = E_{3 \cdot 1}^{3} E_{3 \cdot 2}^{1} - E_{3 \cdot 1}^{1} E_{3 \cdot 2}^{3} = 0 * 1 - 0 = 0$$
$$E_{3 \cdot 1}^{3} = E_{3 \cdot 2}^{1} E_{3 \cdot 3}^{2} - E_{3 \cdot 2}^{2} E_{3 \cdot 3}^{1} = 1 * 0 - 0 = 0$$
$$E_{3 \cdot 2}^{3} = E_{3 \cdot 3}^{1} E_{3 \cdot 1}^{2} - E_{3 \cdot 3}^{2} E_{3 \cdot 1}^{1} = 0 * 1 - 0 * 1 = 0$$
$$E_{3 \cdot 3}^{3} = E_{3 \cdot 1}^{1} E_{3 \cdot 2}^{2} - E_{3 \cdot 1}^{2} E_{3 \cdot 2}^{1} = 1 * 0 - 1 * 1 = -1$$

$\square$

*Замечание* 8.2. Мы можем убедиться непосредственной проверкой, что $\overline{E}_3$ - линейный автоморфизм. Пусть

$$a = a^0 + a^1 i + a^2 j + a^3 k$$
$$b = b^0 + b^1 i + b^2 j + b^3 k$$

Тогда

$$ab = a^0 b^0 - a^1 b^1 - a^2 b^2 - a^3 b^3 + (a^0 b^1 + a^1 b^0 + a^2 b^3 - a^3 b^2)i$$
$$+ (a^0 b^2 + a^2 b^0 + a^3 b^1 - a^1 b^3)j + (a^0 b^3 + a^3 b^0 + a^1 b^2 - a^2 b^1)k$$
$$\overline{E}_3 \circ a = a^0 + a^2 i + a^1 j - a^3 k$$
$$\overline{E}_3 \circ b = b^0 + b^2 i + b^1 j - b^3 k$$
$$(\overline{E}_3 \circ a)(\overline{E}_3 \circ b) = a^0 b^0 - a^2 b^2 - a^1 b^1 - a^3 b^3 + (a^0 b^2 + a^2 b^0 - a^1 b^3 + a^3 b^1)i$$
$$+ (a^0 b^1 + a^1 b^0 - a^3 b^2 + a^2 b^3)j + (-a^0 b^3 - a^3 b^0 + a^2 b^1 - a^1 b^2)k$$
$$= \overline{E}_3 \circ (ab)$$

□

**Теорема 8.3.** *Отображение $\overline{E}_3$ имеет вид*

$$\overline{E}_3 \circ a = -\frac{1}{2}(iai + jaj + iaj + jai)$$

*Доказательство.* Согласно теореме [3]-3.3.4, стандартные компоненты отображения $\overline{E}_3$ имеют вид

$$E_1^{11} = -\tfrac{1}{2} \quad E_1^{22} = -\tfrac{1}{2} \quad E_1^{21} = -\tfrac{1}{2} \quad E_1^{12} = -\tfrac{1}{2}$$

□

**Теорема 8.4.** *Мы можем отождествить отображение*

$$a \circ \overline{E}_3 : H \to H \quad a \in H$$

*и матрицу*

(8.1) $$J_{3l \cdot a} = \begin{pmatrix} a^0 & -a^2 & -a^1 & a^3 \\ a^1 & -a^3 & a^0 & -a^2 \\ a^2 & a^0 & a^3 & a^1 \\ a^3 & a^1 & -a^2 & -a^0 \end{pmatrix}$$

(8.2) $$J_{3l \cdot a} = J_{l \cdot a} \circ E_3$$

*Доказательство.* Произведение кватернионов

$$a = a^0 + a^1 i + a^2 j + a^3 k$$

и

$$\overline{E}_3 \circ x = x^0 + x^2 i + x^1 j - x^3 k$$

имеет вид

$$a \circ \overline{E}_3 \circ x = a^0 x^0 - a^1 x^2 - a^2 x^1 + a^3 x^3 + (a^0 x^2 + a^1 x^0 - a^2 x^3 - a^3 x^1)i$$
$$+ (a^0 x^1 + a^2 x^0 + a^3 x^2 + a^1 x^3)j + (-a^0 x^3 + a^3 x^0 + a^1 x^1 - a^2 x^2)k$$

Следовательно, отображение $f_a \circ x = a \circ \overline{E}_1 \circ x$ имеет матрицу Якоби (8.1). Равенство (8.2) следует из цепочки равенств

$$J_{l \cdot a \circ} {}^\circ E_3 = \begin{pmatrix} a^0 & -a^1 & -a^2 & -a^3 \\ a^1 & a^0 & -a^3 & a^2 \\ a^2 & a^3 & a^0 & -a^1 \\ a^3 & -a^2 & a^1 & a^0 \end{pmatrix} \circ \begin{pmatrix} 1 & 0 & 0 & 0 \\ 0 & 0 & 1 & 0 \\ 0 & 1 & 0 & 0 \\ 0 & 0 & 0 & -1 \end{pmatrix}$$

$$= \begin{pmatrix} a^0 & -a^2 & -a^1 & a^3 \\ a^1 & -a^3 & a^0 & -a^2 \\ a^2 & a^0 & a^3 & a^1 \\ a^3 & a^1 & -a^2 & -a^0 \end{pmatrix} = J_{3l \cdot a}$$

□

## 9. Отображение $\overline{I}$

Антилинейный автоморфизм

$$\overline{I} : H \to H \quad \overline{I} \circ x = x^*$$

$$I^0_0 = 1 \quad I^1_1 = -1 \quad I^2_2 = -1 \quad I^3_3 = -1$$

$$I = \begin{pmatrix} 1 & 0 & 0 & 0 \\ 0 & -1 & 0 & 0 \\ 0 & 0 & -1 & 0 \\ 0 & 0 & 0 & -1 \end{pmatrix}$$

называется сопряжением алгебры кватернионов.

**Теорема 9.1.** *Мы можем отождествить отображение*

$$a \circ \overline{I} : H \to H \quad a \in H$$

*и матрицу*

(9.1) $$I_{l \cdot a} = \begin{pmatrix} a^0 & a^1 & a^2 & a^3 \\ a^1 & -a^0 & a^3 & -a^2 \\ a^2 & -a^3 & -a^0 & a^1 \\ a^3 & a^2 & -a^1 & -a^0 \end{pmatrix}$$

(9.2) $$I_{l \cdot a} = J_{l \cdot a \circ} {}^\circ I$$

*Доказательство.* Произведение кватернионов

$$a = a^0 + a^1 i + a^2 j + a^3 k$$

и

$$I \circ x = x^0 - x^1 i - x^2 j - x^3 k$$

имеет вид

$$a \circ I \circ x = a^0 x^0 + a^1 x^1 + a^2 x^2 + a^3 x^3 + (-a^0 x^1 + a^1 x^0 - a^2 x^3 + a^3 x^2)i$$
$$+ (-a^0 x^2 + a^2 x^0 - a^3 x^1 + a^1 x^3)j + (-a^0 x^3 + a^3 x^0 - a^1 x^2 + a^2 x^1)k$$

Следовательно, отображение $f_a \circ x = a \circ I \circ x$ имеет матрицу Якоби (9.1). Равенство (9.2) следует из цепочки равенств

$$J_{l \cdot a \circ}{}^{\circ} I = \begin{pmatrix} a^0 & -a^1 & -a^2 & -a^3 \\ a^1 & a^0 & -a^3 & a^2 \\ a^2 & a^3 & a^0 & -a^1 \\ a^3 & -a^2 & a^1 & a^0 \end{pmatrix} \circ {}^{\circ} \begin{pmatrix} 1 & 0 & 0 & 0 \\ 0 & -1 & 0 & 0 \\ 0 & 0 & -1 & 0 \\ 0 & 0 & 0 & -1 \end{pmatrix}$$

$$= \begin{pmatrix} a^0 & a^1 & a^2 & a^3 \\ a^1 & -a^0 & a^3 & -a^2 \\ a^2 & -a^3 & -a^0 & a^1 \\ a^3 & a^2 & -a^1 & -a^0 \end{pmatrix} = I_{l \cdot a}$$

□

**Теорема 9.2.** *Мы можем отождествить отображение*

$$a \star \overline{I} : H \to H \quad a \in H$$

*и матрицу*

(9.3) $$I_{r \cdot a} = \begin{pmatrix} a^0 & a^1 & a^2 & a^3 \\ a^1 & -a^0 & -a^3 & a^2 \\ a^2 & a^3 & -a^0 & -a^1 \\ a^3 & -a^2 & a^1 & -a^0 \end{pmatrix}$$

(9.4) $$I_{r \cdot a} = J_{r \cdot a \circ}{}^{\circ} I$$

*Доказательство.* Произведение кватернионов

$$I \circ x = x^0 - x^1 i - x^2 j - x^3 k$$

и

$$a = a^0 + a^1 i + a^2 j + a^3 k$$

имеет вид

$$(I \circ x)a = x^0 a^0 + x^1 a^1 + x^2 a^2 + x^3 a^3 + (x^0 a^1 - x^1 a^0 - x^2 a^3 + x^3 a^2)i$$
$$+ (x^0 a^2 - x^2 a^0 - x^3 a^1 + x^1 a^3)j + (x^0 a^3 - x^3 a^0 - x^1 a^2 + x^2 a^1)k$$

Следовательно, отображение $f_a \circ x = a \circ I \circ x$ имеет матрицу Якоби (9.3). Равенство (9.4) следует из цепочки равенств

$$J_{r \cdot a \circ}{}^\circ I = \begin{pmatrix} a^0 & -a^1 & -a^2 & -a^3 \\ a^1 & a^0 & a^3 & -a^2 \\ a^2 & -a^3 & a^0 & a^1 \\ a^3 & a^2 & -a^1 & a^0 \end{pmatrix} \circ^\circ \begin{pmatrix} 1 & 0 & 0 & 0 \\ 0 & -1 & 0 & 0 \\ 0 & 0 & -1 & 0 \\ 0 & 0 & 0 & -1 \end{pmatrix}$$

$$= \begin{pmatrix} a^0 & a^1 & a^2 & a^3 \\ a^1 & -a^0 & -a^3 & a^2 \\ a^2 & a^3 & -a^0 & -a^1 \\ a^3 & -a^2 & a^1 & -a^0 \end{pmatrix} = I_{r \cdot a}$$

$\square$

## 10. Отображение $\overline{I}_1$

**Теорема 10.1.** *Если отображение*

$$f_2 : H \to H$$

*является антилинейным автоморфизмом алгебры кватернионов и отображение*

$$f_1 : H \to H$$

*является линейным автоморфизмом алгебры кватернионов, то отображение $f_2 \circ f_1$ является антилинейным автоморфизмом алгебры кватернионов.*

*Доказательство.* Так как алгебра кватернионов является $R$-алгеброй, $R \subset \operatorname{Re} H$, то, для $r \in R$, $r^* = r$. Следовательно,

$$f_2 \circ f_1 \circ (ra) = f_2 \circ (rf_1 \circ a) = (f_2 \circ f_1 \circ a)r^*$$

Утверждение теоремы является следствием следующих равенств

$$f_2 \circ f_1 \circ (a+b) = f_2 \circ (f_1 \circ a + f_1 \circ b) = f_2 \circ f_1 \circ a + f_2 \circ f_1 \circ b$$
$$f_2 \circ f_1 \circ (ab) = f_2((f_1 \circ a)(f_1 \circ b)) = (f_2 \circ f_1 \circ b)(f_2 \circ f_1 \circ a)$$

$\square$

**Теорема 10.2.** *Отображение*

$$\overline{I}_1 : H \to H \quad \overline{I}_1 \circ x = x^0 - x^2 \overline{e}_1 - x^3 \overline{e}_2 - x^1 \overline{e}_3$$

$$I_{1 \cdot 0}^{\ 0} = 1 \quad I_{1 \cdot 2}^{\ 1} = -1 \quad I_{1 \cdot 3}^{\ 2} = -1 \quad I_{1 \cdot 1}^{\ 3} = -1$$

$$I_1 = \begin{pmatrix} 1 & 0 & 0 & 0 \\ 0 & 0 & -1 & 0 \\ 0 & 0 & 0 & -1 \\ 0 & -1 & 0 & 0 \end{pmatrix}$$

*является антилинейным автоморфизмом алгебры кватернионов.*

*Доказательство.* Утверждение теоремы является следствием теоремы 10.1 и равенства
$$\overline{I}_1 = \overline{I} \circ \overline{E}_1$$
□

**Теорема 10.3.** *Мы можем отождествить отображение*
$$a \circ \overline{I}_1 : H \to H \quad a \in H$$
*и матрицу*

(10.1) $$I_{1l \cdot a} = \begin{pmatrix} a^0 & a^3 & a^1 & a^2 \\ a^1 & -a^2 & -a^0 & a^3 \\ a^2 & a^1 & -a^3 & -a^0 \\ a^3 & -a^0 & a^2 & -a^1 \end{pmatrix}$$

(10.2) $$I_{1l \cdot a} = J_{1l \cdot a \circ} {}^\circ I_1$$

*Доказательство.* Произведение кватернионов
$$a = a^0 + a^1 i + a^2 j + a^3 k$$
и
$$\overline{I}_1 \circ x = x^0 - x^2 i - x^3 j - x^1 k$$
имеет вид
$$a \circ \overline{I}_1 \circ x = a^0 x^0 + a^1 x^2 + a^2 x^3 + a^3 x^1 + (-a^0 x^2 + a^1 x^0 - a^2 x^1 + a^3 x^3)i$$
$$+ (-a^0 x^3 + a^2 x^0 - a^3 x^2 + a^1 x^1)j + (-a^0 x^1 + a^3 x^0 - a^1 x^3 + a^2 x^2)k$$

Следовательно, отображение $f_a \circ x = a \circ \overline{I}_1 \circ x$ имеет матрицу Якоби (10.1). Равенство (10.2) следует из цепочки равенств

$$J_{l \cdot a \circ} {}^\circ I_1 = \begin{pmatrix} a^0 & -a^1 & -a^2 & -a^3 \\ a^1 & a^0 & -a^3 & a^2 \\ a^2 & a^3 & a^0 & -a^1 \\ a^3 & -a^2 & a^1 & a^0 \end{pmatrix} {}_\circ^\circ \begin{pmatrix} 1 & 0 & 0 & 0 \\ 0 & 0 & -1 & 0 \\ 0 & 0 & 0 & -1 \\ 0 & -1 & 0 & 0 \end{pmatrix}$$

$$= \begin{pmatrix} a^0 & a^3 & a^1 & a^2 \\ a^1 & -a^2 & -a^0 & a^3 \\ a^2 & a^1 & -a^3 & -a^0 \\ a^3 & -a^0 & a^2 & -a^1 \end{pmatrix} = I_{1l \cdot a}$$

□

**Теорема 10.4.** *Мы можем отождествить отображение*
$$a \star \overline{I}_1 : H \to H \quad a \in H$$

и матрицу

(10.3) $$I_{1r\cdot a} = \begin{pmatrix} a^0 & a^3 & a^1 & a^2 \\ a^1 & a^2 & -a^0 & -a^3 \\ a^2 & -a^1 & a^3 & -a^0 \\ a^3 & -a^0 & -a^2 & a^1 \end{pmatrix}$$

(10.4) $$I_{1r\cdot a} = J_{1r\cdot a\circ}{}^\circ I_1$$

*Доказательство.* Произведение кватернионов
$$\overline{I}_1 \circ x = x^0 - x^2 i - x^3 j - x^1 k$$
и
$$a = a^0 + a^1 i + a^2 j + a^3 k$$
имеет вид
$$(I \circ x)a = x^0 a^0 + x^2 a^1 + x^3 a^2 + x^1 a^3 + (x^0 a^1 - x^2 a^0 - x^3 a^3 + x^1 a^2)i$$
$$+ (x^0 a^2 - x^3 a^0 - x^1 a^1 + x^2 a^3)j + (x^0 a^3 - x^1 a^0 - x^2 a^2 + x^3 a^1)k$$

Следовательно, отображение $f_a \circ x = a \star \overline{I}_1 \circ x$ имеет матрицу Якоби (10.3). Равенство (10.4) следует из цепочки равенств

$$J_{r\cdot a\circ}{}^\circ I_1 = \begin{pmatrix} a^0 & -a^1 & -a^2 & -a^3 \\ a^1 & a^0 & a^3 & -a^2 \\ a^2 & -a^3 & a^0 & a^1 \\ a^3 & a^2 & -a^1 & a^0 \end{pmatrix} {}_\circ{}^\circ \begin{pmatrix} 1 & 0 & 0 & 0 \\ 0 & 0 & -1 & 0 \\ 0 & 0 & 0 & -1 \\ 0 & -1 & 0 & 0 \end{pmatrix}$$
$$= \begin{pmatrix} a^0 & a^3 & a^1 & a^2 \\ a^1 & a^2 & -a^0 & -a^3 \\ a^2 & -a^1 & a^3 & -a^0 \\ a^3 & -a^0 & -a^2 & a^1 \end{pmatrix} = I_{1r\cdot a}$$

□

## 11. Отображение $\overline{I}_2$

**Теорема 11.1.** *Отображение*
$$\overline{I}_2 : H \to H \quad \overline{I}_2 \circ x = x^0 - x^3 \overline{e}_1 - x^1 \overline{e}_2 - x^2 \overline{e}_3$$
$$I_{2\cdot 0}^{\,0} = 1 \quad I_{2\cdot 3}^{\,1} = -1 \quad I_{2\cdot 1}^{\,2} = -1 \quad I_{2\cdot 2}^{\,3} = -1$$
$$I_2 = \begin{pmatrix} 1 & 0 & 0 & 0 \\ 0 & 0 & 0 & -1 \\ 0 & -1 & 0 & 0 \\ 0 & 0 & -1 & 0 \end{pmatrix}$$
*является антилинейным автоморфизмом алгебры кватернионов.*

*Доказательство.* Утверждение теоремы является следствием теоремы 10.1 и равенства
$$\overline{I}_2 = \overline{I} \circ \overline{E}_2$$

□

**Теорема 11.2.** *Мы можем отождествить отображение*
$$a \circ \overline{I}_2 : H \to H \quad a \in H$$
*и матрицу*

(11.1)
$$I_{2l \cdot a} = \begin{pmatrix} a^0 & a^2 & a^3 & a^1 \\ a^1 & a^3 & -a^2 & -a^0 \\ a^2 & -a^0 & a^1 & -a^3 \\ a^3 & -a^1 & -a^0 & a^2 \end{pmatrix}$$

(11.2)
$$I_{2l \cdot a} = J_{l \cdot a} \circ° I_2$$

*Доказательство.* Произведение кватернионов
$$a = a^0 + a^1 i + a^2 j + a^3 k$$
и
$$\overline{I}_2 \circ x = x^0 - x^3 i - x^1 j - x^2 k$$
имеет вид
$$a \circ \overline{I}_2 \circ x = a^0 x^0 + a^1 x^3 + a^2 x^1 + a^3 x^2 + (-a^0 x^3 + a^1 x^0 - a^2 x^2 + a^3 x^1)i$$
$$+ (-a^0 x^1 + a^2 x^0 - a^3 x^3 + a^1 x^2)j + (-a^0 x^2 + a^3 x^0 - a^1 x^1 + a^2 x^3)k$$

Следовательно, отображение $f_a \circ x = a \circ \overline{I}_2 \circ x$ имеет матрицу Якоби (11.1). Равенство (11.2) следует из цепочки равенств

$$J_{l \cdot a} \circ° I_2 = \begin{pmatrix} a^0 & -a^1 & -a^2 & -a^3 \\ a^1 & a^0 & -a^3 & a^2 \\ a^2 & a^3 & a^0 & -a^1 \\ a^3 & -a^2 & a^1 & a^0 \end{pmatrix} \circ° \begin{pmatrix} 1 & 0 & 0 & 0 \\ 0 & 0 & 0 & -1 \\ 0 & -1 & 0 & 0 \\ 0 & 0 & -1 & 0 \end{pmatrix}$$
$$= \begin{pmatrix} a^0 & a^2 & a^3 & a^1 \\ a^1 & a^3 & -a^2 & -a^0 \\ a^2 & -a^0 & a^1 & -a^3 \\ a^3 & -a^1 & -a^0 & a^2 \end{pmatrix} = I_{2l \cdot a}$$

□

## 12. Линейное отображение

**Теорема 12.1.** *Линейное отображение алгебры кватернионов*
$$\overline{f} : H \to H$$
*имеет единственное разложение*

(12.1) $$\overline{f} = a_0 \circ \overline{E} + a_1 \circ \overline{E}_1 + a_2 \circ \overline{E}_3 + a_3 \circ \overline{I}$$

*При этом*

(12.2) $$a_0^0 = \frac{1}{2}(f_0^0 + f_1^1 - f_2^1 + f_3^0)$$

(12.3) $$a_0^1 = \frac{1}{2}(f_0^1 - f_1^0 + f_2^0 + f_3^1)$$

(12.4) $$a_0^2 = \frac{1}{2}(f_0^2 - f_1^3 + f_2^3 + f_3^2)$$

(12.5) $$a_0^3 = \frac{1}{2}(f_0^3 + f_1^2 - f_2^2 + f_3^3)$$

(12.6) $$a_1^0 = \frac{1}{2}(-f_1^0 - f_1^1 + f_2^0 + f_2^1 - f_3^0 + f_3^1 + f_3^2 + f_3^3)$$

(12.7) $$a_1^1 = \frac{1}{2}(f_1^0 - f_1^1 - f_2^0 + f_2^1 - f_3^0 - f_3^1 - f_3^2 + f_3^3)$$

(12.8) $$a_1^2 = \frac{1}{2}(-f_1^2 + f_1^3 + f_2^2 - f_2^3 - f_3^0 + f_3^1 - f_3^2 - f_3^3)$$

(12.9) $$a_1^3 = \frac{1}{2}(-f_1^2 - f_1^3 + f_2^2 + f_2^3 - f_3^0 - f_3^1 + f_3^2 - f_3^3)$$

(12.10) $$a_2^0 = \frac{1}{2}(f_1^0 + f_1^2 - f_2^0 + f_2^1 - f_3^1 - f_3^2)$$

(12.11) $$a_2^1 = \frac{1}{2}(f_1^1 + f_1^3 - f_2^0 - f_2^1 + f_3^0 - f_3^3)$$

(12.12) $$a_2^2 = \frac{1}{2}(-f_1^0 + f_1^2 - f_2^2 - f_2^3 + f_3^0 + f_3^3)$$

(12.13) $$a_2^3 = \frac{1}{2}(-f_1^1 + f_1^3 + f_2^2 - f_2^3 + f_3^1 - f_3^2)$$

(12.14) $$a_3^0 = \frac{1}{2}(f_0^0 - f_1^2 - f_2^1 - f_3^3)$$

(12.15) $$a_3^1 = \frac{1}{2}(f_0^1 - f_1^3 + f_2^0 + f_3^2)$$

(12.16) $$a_3^2 = \frac{1}{2}(f_0^2 + f_1^0 + f_2^3 - f_3^1)$$

(12.17) $$a_3^3 = \frac{1}{2}(f_0^3 + f_1^1 - f_2^2 + f_3^0)$$

*Доказательство.* Линейное отображение (12.1) имеет матрицу

$$
(12.18) \quad
\begin{pmatrix}
a_0^0 & -a_0^1 & -a_0^2 & -a_0^3 \\
a_0^1 & a_0^0 & -a_0^3 & a_0^2 \\
a_0^2 & a_0^3 & a_0^0 & -a_0^1 \\
a_0^3 & -a_0^2 & a_0^1 & a_0^0
\end{pmatrix}
+
\begin{pmatrix}
a_1^0 & -a_1^3 & -a_1^1 & -a_1^2 \\
a_1^1 & a_1^2 & a_1^0 & -a_1^3 \\
a_1^2 & -a_1^1 & a_1^3 & a_1^0 \\
a_1^3 & a_1^0 & -a_1^2 & a_1^1
\end{pmatrix}
$$
$$
+
\begin{pmatrix}
a_2^0 & -a_2^2 & -a_2^1 & a_2^3 \\
a_2^1 & -a_2^3 & a_2^0 & -a_2^2 \\
a_2^2 & a_2^0 & a_2^3 & a_2^1 \\
a_2^3 & a_2^1 & -a_2^2 & -a_2^0
\end{pmatrix}
+
\begin{pmatrix}
a_3^0 & a_3^1 & a_3^2 & a_3^3 \\
a_3^1 & -a_3^0 & a_3^3 & -a_3^2 \\
a_3^2 & -a_3^3 & -a_3^0 & a_3^1 \\
a_3^3 & a_3^2 & -a_3^1 & -a_3^0
\end{pmatrix}
$$

Из сравнения матрицы отображения $\overline{f}$ и матрицы (12.18), мы получим систему линейных уравнений

(12.19) $\quad f_0^0 = a_0^0 + a_1^0 + a_2^0 + a_3^0$ $\qquad$ (12.27) $\quad f_0^2 = a_0^2 + a_1^2 + a_2^2 + a_3^2$

(12.20) $\quad f_1^0 = -a_0^1 - a_1^3 - a_2^2 + a_3^1$ $\qquad$ (12.28) $\quad f_1^2 = a_0^3 - a_1^1 + a_2^0 - a_3^3$

(12.21) $\quad f_2^0 = -a_0^2 - a_1^1 - a_2^1 + a_3^2$ $\qquad$ (12.29) $\quad f_2^2 = a_0^0 + a_1^3 + a_2^3 - a_3^0$

(12.22) $\quad f_3^0 = -a_0^3 - a_1^2 + a_2^3 + a_3^3$ $\qquad$ (12.30) $\quad f_3^2 = -a_0^1 + a_1^0 + a_2^1 + a_3^1$

(12.23) $\quad f_0^1 = a_0^0 + a_1^1 + a_2^1 + a_3^1$ $\qquad$ (12.31) $\quad f_0^3 = a_0^3 + a_1^3 + a_2^3 + a_3^3$

(12.24) $\quad f_1^1 = a_0^0 + a_1^2 - a_2^3 - a_3^0$ $\qquad$ (12.32) $\quad f_1^3 = -a_0^2 + a_1^0 + a_2^1 + a_3^2$

(12.25) $\quad f_2^1 = -a_0^3 + a_1^0 + a_2^0 + a_3^3$ $\qquad$ (12.33) $\quad f_2^3 = a_0^1 - a_1^2 - a_2^2 - a_3^1$

(12.26) $\quad f_3^1 = a_0^2 - a_1^3 - a_2^2 - a_3^2$ $\qquad$ (12.34) $\quad f_3^3 = a_0^0 + a_1^1 - a_2^2 - a_3^0$

Для решения этой системы линейный уравнений я написал программу на C#. Решение легко проверить непосредственной подстановкой. $\square$

Из теоремы [3]-2.6.4 следует, что множество линейных эндоморфизмов $\mathcal{L}(H;H)$ алгебры кватернионов $H$ изоморфно тензорному произведению $H \otimes H$. Теорема 12.1 утверждает, что мы можем рассматривать модуль $\mathcal{L}(H;H)$ как $H\star$-векторное пространство с базисом

$$(\overline{E}, \overline{E}_1, \overline{E}_3, \overline{I})$$

**Пример 12.2.** Согласно теореме 7.4, отображение

$$a \circ \overline{E}_2 : H \to H \quad a \in H$$

имеет матрицу

$$
\begin{pmatrix}
a^0 & -a^2 & -a^3 & -a^1 \\
a^1 & -a^3 & a^2 & a^0 \\
a^2 & a^0 & -a^1 & a^3 \\
a^3 & a^1 & a^0 & -a^2
\end{pmatrix}
$$

Согласно теореме 12.1,

$$a_0^0 = \frac{1}{2}(a^0 + (-a^3) - a^2 + (-a^1))$$

$$a_0^1 = \frac{1}{2}(a^1 - (-a^2) + (-a^3) + a^0)$$

$$a_0^2 = \frac{1}{2}(a^2 - a^1 + a^0 + a^3)$$

$$a_0^3 = \frac{1}{2}(a^3 + a^0 - (-a^1) + (-a^2))$$

$$a_1^0 = \frac{1}{2}(-(-a^2) - (-a^3) + (-a^3) + a^2 - (-a^1) + a^0 + a^3 + (-a^2))$$

$$a_1^1 = \frac{1}{2}((-a^2) - (-a^3) - (-a^3) + a^2 - (-a^1) - a^0 - a^3 + (-a^2))$$

$$a_1^2 = \frac{1}{2}(-a^0 + a^1 + (-a^1) - a^0 - (-a^1) + a^0 - a^3 - (-a^2))$$

$$a_1^3 = \frac{1}{2}(-a^0 - a^1 + (-a^1) + a^0 - (-a^1) - a^0 + a^3 - (-a^2))$$

$$a_2^0 = \frac{1}{2}((-a^2) + a^0 - (-a^3) + a^2 - a^0 - a^3)$$

$$a_2^1 = \frac{1}{2}((-a^3) + a^1 - (-a^3) - a^2 + (-a^1) - (-a^2))$$

$$a_2^2 = \frac{1}{2}(-(-a^2) + a^0 - (-a^1) - a^0 + (-a^1) + (-a^2))$$

$$a_2^3 = \frac{1}{2}(-(-a^3) + a^1 + (-a^1) - a^0 + a^0 - a^3)$$

$$a_3^0 = \frac{1}{2}(a^0 - a^0 - a^2 - (-a^2))$$

$$a_3^1 = \frac{1}{2}(a^1 - a^1 + (-a^3) + a^3)$$

$$a_3^2 = \frac{1}{2}(a^2 + (-a^2) + a^0 - a^0)$$

$$a_3^3 = \frac{1}{2}(a^3 + (-a^3) - (-a^1) + (-a^1))$$

Следовательно, $a_2 = a_3 = 0$. $\square$

**Пример 12.3.** Согласно теореме 5.2, отображение

$$a \star \overline{E} : H \to H \quad a \in H$$

имеет матрицу

$$\begin{pmatrix} a^0 & -a^1 & -a^2 & -a^3 \\ a^1 & a^0 & a^3 & -a^2 \\ a^2 & -a^3 & a^0 & a^1 \\ a^3 & a^2 & -a^1 & a^0 \end{pmatrix}$$

Согласно теореме 12.1,

$$a_0^0 = \frac{1}{2}(a^0 + a^0 - a^3 + (-a^3))$$
$$a_0^1 = \frac{1}{2}(a^1 - (-a^1) + (-a^2) + (-a^2))$$
$$a_0^2 = \frac{1}{2}(a^2 - a^2 + (-a^1) + a^1)$$
$$a_0^3 = \frac{1}{2}(a^3 + (-a^3) - a^0 + a^0)$$

$$a_1^0 = \frac{1}{2}(-(-a^1) - a^0 + (-a^2) + a^3 - (-a^3) + (-a^2) + a^1 + a^0)$$
$$a_1^1 = \frac{1}{2}((-a^1) - a^0 - (-a^2) + a^3 - (-a^3) - (-a^2) - a^1 + a^0)$$
$$a_1^2 = \frac{1}{2}(-(-a^3) + a^2 + a^0 - (-a^1) - (-a^3) + (-a^2) - a^1 - a^0)$$
$$a_1^3 = \frac{1}{2}(-(-a^3) - a^2 + a^0 + (-a^1) - (-a^3) - (-a^2) + a^1 - a^0)$$

$$a_2^0 = \frac{1}{2}((-a^1) + (-a^3) - (-a^2) + a^3 - (-a^2) - a^1)$$
$$a_2^1 = \frac{1}{2}(a^0 + a^2 - (-a^2) - a^3 + (-a^3) - a^0)$$
$$a_2^2 = \frac{1}{2}(-(-a^1) + (-a^3) - a^0 - (-a^1) + (-a^3) + a^0)$$
$$a_2^3 = \frac{1}{2}(-a^0 + a^2 + a^0 - (-a^1) + (-a^2) - a^1)$$

$$a_3^0 = \frac{1}{2}(a^0 - (-a^3) - a^3 - a^0)$$
$$a_3^1 = \frac{1}{2}(a^1 - a^2 + (-a^2) + a^1)$$
$$a_3^2 = \frac{1}{2}(a^2 + (-a^1) + (-a^1) - (-a^2))$$
$$a_3^3 = \frac{1}{2}(a^3 + a^0 - a^0 + (-a^3))$$

Это очень важно, $a_3 \ne 0$. $\square$

**Пример 12.4.** Согласно теоремам 5.1, 5.2, отображение

$$f : H \to H \quad f = a \circ E + a \star E \quad a \in H$$
$$f \circ x = ax + xa$$

имеет матрицу

$$\begin{pmatrix} a^0 & -a^1 & -a^2 & -a^3 \\ a^1 & a^0 & -a^3 & a^2 \\ a^2 & a^3 & a^0 & -a^1 \\ a^3 & -a^2 & a^1 & a^0 \end{pmatrix} + \begin{pmatrix} a^0 & -a^1 & -a^2 & -a^3 \\ a^1 & a^0 & a^3 & -a^2 \\ a^2 & -a^3 & a^0 & a^1 \\ a^3 & a^2 & -a^1 & a^0 \end{pmatrix}$$

$$= \begin{pmatrix} 2a^0 & -2a^1 & -2a^2 & -2a^3 \\ 2a^1 & 2a^0 & 0 & 0 \\ 2a^2 & 0 & 2a^0 & 0 \\ 2a^3 & 0 & 0 & 2a^0 \end{pmatrix}$$

Согласно теореме 12.1,

$$a_0^0 = \frac{1}{2}(2a^0 + 2a^0 - 0 + (-2a^3))$$
$$a_0^1 = \frac{1}{2}(2a^1 - (-2a^1) + (-2a^2) + 0)$$
$$a_0^2 = \frac{1}{2}(2a^2 - 0 + 0 + 0)$$
$$a_0^3 = \frac{1}{2}(2a^3 + 0 - 2a^0 + 2a^0)$$

$$a_1^0 = \frac{1}{2}(-(-2a^1) - 2a^0 + (-2a^2) + 0 - (-2a^3) + 0 + 0 + 2a^0)$$
$$a_1^1 = \frac{1}{2}((-2a^1) - 2a^0 - (-2a^2) + 0 - (-2a^3) - 0 - 0 + 2a^0)$$
$$a_1^2 = \frac{1}{2}(-0 + 0 + 2a^0 - 0 - (-2a^3) + 0 - 0 - 2a^0)$$
$$a_1^3 = \frac{1}{2}(-0 - 0 + 2a^0 + 0 - (-2a^3) - 0 + 0 - 2a^0)$$

$$a_2^0 = \frac{1}{2}((-2a^1) + 0 - (-2a^2) + 0 - 0 - 0)$$
$$a_2^1 = \frac{1}{2}(2a^0 + 0 - (-2a^2) - 0 + (-2a^3) - 2a^0)$$
$$a_2^2 = \frac{1}{2}(-(-2a^1) + 0 - 2a^0 - 0 + (-2a^3) + 2a^0)$$
$$a_2^3 = \frac{1}{2}(-2a^0 + 0 + 2a^0 - 0 + 0 - 0)$$

$$a_3^0 = \frac{1}{2}(2a^0 - 0 - 0 - 2a^0)$$
$$a_3^1 = \frac{1}{2}(2a^1 - 0 + (-2a^2) + 0)$$
$$a_3^2 = \frac{1}{2}(2a^2 + (-2a^1) + 0 - 0)$$
$$a_3^3 = \frac{1}{2}(2a^3 + 2a^0 - 2a^0 + (-2a^3))$$

Это очень важно, $a_3 \neq 0$. □

## 13. Список литературы

# 14. Специальные символы и обозначения